\documentclass[hidelinks,12pt]{amsart}
\usepackage{bbm}
\usepackage{libertine}
\usepackage{enumitem}
\usepackage[T1]{fontenc}
\usepackage{footmisc}
\usepackage{stmaryrd}
\usepackage{pdfpages}
\usepackage{color}
\DeclareFontFamily{OT1}{pzc}{}
\DeclareFontShape{OT1}{pzc}{m}{it}{<-> s * [1.100] pzcmi7t}{}
\DeclareMathAlphabet{\mathpzc}{OT1}{pzc}{m}{it}\DeclareMathAlphabet{\mathcal}{OMS}{cmsy}{m}{n}\usepackage{fullpage}
\usepackage{mdframed}
\usepackage{graphicx}
\usepackage{calligra}
\usepackage{wasysym}
\usepackage{csquotes}
\usepackage{bm}
\usepackage{amsmath}
\usepackage[all]{xy}
\usepackage{amssymb}
\usepackage{amsthm}
\usepackage{mathrsfs}
\usepackage[OT2,T1]{fontenc}
\usepackage{centernot}
\usepackage{hyperref}

\newcommand{\DMO}[2]{\DeclareMathOperator{#1}{#2}}
\DMO{\BK}{BK}
\DMO{\FL}{FL}
\DMO{\Ann}{Ann}
\DMO{\std}{std}
\DMO{\antidiag}{antidiag}
\DMO{\locadm}{loc.adm}
\DMO{\Inj}{Inj}
\DMO{\LL}{LL}
\DMO{\Dmod}{\emph{D }-mod}
\DMO{\univ}{univ}
\DMO{\Fitt}{Fitt}
\DMO{\WD}{WD}
\DMO{\geom}{geom}
\DMO{\Fl}{Fl}
\DMO{\grad}{grad}
\DMO{\labmda}{\lambda}
\DMO{\Iw}{Iw}
\DMO{\tor}{tor}
\DMO{\coh}{coh}
\DMO{\vol}{vol}
\DMO{\semsim}{ss}
\DMO{\free}{free}
\DMO{\Alg}{Alg}
\DMO{\oth}{otherwise}
\DMO{\Ber}{Ber}
\DMO{\Diff}{Diff}
\DMO{\br}{br}
\DMO{\Isot}{Isot}
\DMO{\prim}{prim}
\DMO{\RAH}{RAH}
\DMO{\Sets}{Sets}
\DMO{\cone}{cone}
\DMO{\Grps}{Grps}
\DMO{\Dec}{Dec}
\DMO{\Flat}{Flat}
\DMO{\AbGps}{AbGps}
\DMO{\Sch}{Sch}
\DMO{\AH}{AH}
\DMO{\cl}{cl}
\DMO{\sk}{sk}
\DMO{\HC}{HC}
\DMO{\cosk}{sk}
\DMO{\ur}{ur}
\DMO{\LocSys}{LocSys}
\DMO{\rk}{rk}
\DMO{\NT}{NT}
\DMO{\cork}{cork}
\DMO{\KS}{KS}
\DMO{\MU}{MU}
\DMO{\der}{der}
\DMO{\Art}{Art}
\DMO{\Proj}{Proj}
\DMO{\End}{End}
\DMO{\Betti}{Betti}
\DMO{\Sym}{Sym}
\DMO{\cInd}{cInd}
\DMO{\GL}{GL}
\DMO{\Gal}{Gal}
\DMO{\Br}{Br}
\DMO{\Der}{Der}
\DMO{\Sp}{Sp}
\DMO{\Tan}{Tan}
\DMO{\Spin}{Spin}
\DMO{\Var}{Var}
\DMO{\Nrd}{Nrd}
\DMO{\cusp}{cusp}
\DMO{\Mat}{Mat}
\DMO{\Isom}{Isom}
\DMO{\Stab}{Stab}
\DMO{\SO}{SO}
\DMO{\Res}{Res}
\DMO{\Lie}{Lie}
\DMO{\SU}{SU}
\DMO{\Ad}{Ad}
\DMO{\ad}{ad}
\DMO{\im}{im}
\DMO{\Frob}{Frob}
\DMO{\Fr}{Fr}
\DMO{\red}{red}
\DMO{\an}{an}
\DMO{\Pic}{Pic}
\DMO{\Tor}{Tor}
\DMO{\Hdg}{Hdg}
\DMO{\id}{id}
\DMO{\pr}{pr}
\DMO{\Mor}{Mor}
\DMO{\Ext}{Ext}
\DMO{\ML}{ML}
\DMO{\PGL}{PGL}
\DMO{\SL}{SL}
\DMO{\GU}{GU}
\DMO{\GSp}{GSp}
\DMO{\GSL}{GSL}
\DMO{\Aff}{Aff}
\DMO{\NS}{NS}
\DMO{\gr}{gr}
\DMO{\Ch}{Ch}
\DMO{\QCoh}{QCoh}
\DMO{\Coh}{Coh}
\DMO{\inv}{inv}
\DMO{\Gr}{Gr}
\DMO{\Bun}{Bun}
\DMO{\Hk}{Hk}
\DMO{\GH}{GH}
\DMO{\HT}{HT}
\DMO{\LT}{LT}
\DMO{\Int}{Int}
\DMO{\UU}{U}
\DMO{\OO}{O}
\DMO{\Loc}{Loc}
\DMO{\Conn}{Conn}
\DMO{\sing}{sing}
\DMO{\si}{si}
\DMO{\Sen}{Sen}
\DMO{\MaxSpec}{MaxSpec}
\DMO{\ran}{ran}
\DMO{\coker}{coker}
\DMO{\DIV}{div}
\DMO{\Cl}{Cl}
\DMO{\Frac}{Frac}
\DMO{\VEC}{Vec}
\DMO{\Weil}{Weil}
\DMO{\SPLIT}{split}
\DMO{\Tr}{Tr}
\DMO{\val}{val}
\DMO{\pv}{p.v.}
\DMO{\disc}{disc}
\DMO{\trdeg}{tr.deg}
\DMO{\rad}{rad}
\DMO{\codim}{codim}
\DMO{\dist}{dist}
\DMO{\length}{length}
\DMO{\diam}{diam}
\DMO{\Supp}{Supp}
\DMO{\Ass}{Ass}
\DMO{\ord}{ord}
\DMO{\RE}{Re}
\DMO{\Sh}{Sh}
\DMO{\IM}{Im}
\DMO{\Tot}{Tot}
\DMO{\Bl}{Bl}
\DMO{\lcm}{lcm}
\DMO{\ann}{ann}
\DMO{\arcsinh}{arcsinh}
\DMO{\CHAR}{char}
\DMO{\MOD}{mod}
\DMO{\BB}{BB}
\DMO{\new}{new}
\DMO{\alg}{alg}
\DMO{\Irr}{Irr}
\DMO{\res}{res}
\DMO{\rank}{rank}
\DMO{\naive}{naive}
\DMO{\tors}{tors}
\DMO{\Perf}{Perf}
\DMO{\Sht}{Sht}
\DMO{\Perv}{Perv}
\DMO{\soc}{soc}
\DMO{\Mod}{Mod}
\DMO{\cyc}{cyc}
\DMO{\SC}{sc}
\DMO{\SP}{sp}
\DMO{\Deck}{Deck}
\DMO{\PSL}{PSL}
\DMO{\Area}{Area}
\DMO{\Cont}{Cont}
\DMO{\sgn}{sgn}
\DMO{\Cat}{Cat}
\DMO{\Cov}{Cov}
\DMO{\rig}{rig}
\DMO{\FSch}{FSch}
\DMO{\Rig}{Rig}
\DMO{\Spv}{Spv}
\DMO{\Spa}{Spa}
\DMO{\trace}{trace}
\DMO{\cont}{cont}
\DMO{\aff}{aff}
\DMO{\cor}{cor}
\DMO{\CH}{CH}
\DMO{\Spec}{Spec}
\DMO{\rec}{rec}
\DMO{\LGC}{LGC}
\DMO{\un}{un}
\DMO{\conj}{conj}
\DMO{\Eval}{Eval}
\DMO{\JH}{JH}
\DMO{\can}{can}
\DMO{\Fss}{Fss}
\DMO{\Speh}{Speh}
\DMO{\Ind}{Ind}
\DMO{\ch}{ch}
\DMO{\nr}{nr}
\DMO{\Swan}{Swan}
\DMO{\St}{St}
\DMO{\Ho}{Ho}
\DMO{\HH}{HH}
\DMO{\trop}{trop}
\DMO{\Jac}{Jac}
\DMO{\vir}{vir}
\DMO{\coll}{coll}
\DMO{\reg}{reg}
\DMO{\dlog}{dlog}
\DMO{\Div}{Div}
\DMO{\ab}{ab}
\DMO{\Tam}{Tam}
\DMO{\Ran}{Ran}
\DMO{\IC}{IC}
\DMO{\Sat}{Sat}
\DMO{\Rat}{Rat}
\DMO{\loc}{loc}
\DMO{\ev}{ev}
\DMO{\st}{st}
\DMO{\pst}{pst}
\DMO{\Fil}{Fil}
\DMO{\cris}{cris}
\DMO{\dR}{dR}
\DMO{\Rep}{Rep}
\DMO{\Sel}{Sel}
\DMO{\spec}{spec}
\DMO{\Spf}{Spf}
\DMO{\JL}{JL}
\DMO{\BGL}{BGL}
\DMO{\Arc}{Arc}
\DMO{\MHS}{MHS}
\DMO{\Nm}{Nm}
\DMO{\holim}{holim}
\DMO{\nInd}{nInd}
\DMO{\sSets}{s\textbf{Sets}}
\DMO{\sArt}{s\textbf{Art}}
\DMO{\BDJ}{BDJ}
\DMO{\GV}{GV}
\DMO{\BM}{BM}
\DMO{\Ord}{Ord}
\DMO{\mult}{mult}
\DMO{\WDRep}{WDRep}
\DMO{\Aut}{Aut}
\DMO{\Hom}{Hom}
\DMO{\sph}{sph}
\DMO{\Def}{Def}
\DMO{\GO}{GO}
\DMO{\diag}{diag}
\DMO{\cond}{cond}
\DMO{\ind}{ind}
\DMO{\irr}{irr}
\DMO{\RHom}{RHom}
\DMO{\sm}{sm}
\DMO{\sss}{ss}
\DMO{\sHom}{sHom}
\DMO{\Tran}{Tran}
\DMO{\Rees}{Rees}
\DMO{\lcv}{lcv} 
\DMO{\SN}{SN}
\DMO{\triv}{triv}
\DMO{\height}{ht}
\DMO{\proj}{proj}
\DMO{\Fun}{Fun}
\DMO{\cts}{cts}
\DMO{\Obj}{Obj}
\DMO{\Sing}{Sing}
\DMO{\Pro}{Pro}
\DMO{\Ig}{Ig}
\DMO{\Ha}{Ha}
\DMO{\BC}{BC}
\DMO{\RZ}{RZ}
\DMO{\supp}{supp}
\DMO{\projdim}{proj.dim}
\DMO{\Zar}{Zar}
\DMO{\Ban}{Ban}
\DMO{\LA}{LA}
\DMO{\ess}{ess}
\DMO{\op}{op}
\DMO{\Func}{Func}
\DMO{\Born}{Born}
\DMO{\Comm}{Comm}
\DMO{\Dr}{Dr}
\DMO{\LC}{LC}
\DMO{\nind}{n-ind}
\DMO{\perf}{perf}
\DMO{\charpoly}{char.poly}
\makeatletter
\def\thmhead@plain#1#2#3{%
  \thmname{#1}\thmnumber{\@ifnotempty{#1}{ }\@upn{#2}}%
  \thmnote{ {\the\thm@notefont#3}}}
\let\thmhead\thmhead@plain
\makeatother

\newtheorem*{thm1*}{Theorem}
\newtheorem*{lemma*}{Lemma}
\newtheorem*{defn1*}{Definition}

\newtheorem{thm2}{Theorem}[section]

\newtheorem{innercustomthm}{Theorem}
\newenvironment{customthm}[1]
  {\renewcommand\theinnercustomthm{#1}\innercustomthm}
  {\endinnercustomthm}
\newtheorem{innercustomconj}{Conjecture}
\newenvironment{customconj}[1]
  {\renewcommand\theinnercustomconj{#1}\innercustomconj}
  {\endinnercustomconj}

\newtheorem*{prop*}{Proposition}
\newtheorem{conj2}[thm2]{Conjecture}

\newtheorem{prop2}[thm2]{Proposition}

\newtheorem*{conj*}{Conjecture}

\theoremstyle{definition}

\newtheorem{defthm2}[thm2]{Theorem-Definition}
\newtheorem{defn2}[thm2]{Definition}

\newtheorem*{defn2*}{Definition}

\newtheorem{homework}{}
\newtheorem*{prb*}{Problem}

\newtheorem*{claim*}{Claim}

\newtheorem{rmk2}[thm2]{Remark}

\newtheorem*{rmk*}{Remark}

\newtheorem{exam2}[thm2]{Example}

\newtheorem{ass2}[thm2]{Assumption}

\newtheoremstyle{theoremdd}
  {6pt}
  {6pt}
  {}
  {0pt}
  {\bfseries}
  {.}
  { }
  {\thmname{#1}\thmnumber{ #2}\textnormal{\thmnote{ #3}}}
  \theoremstyle{theoremdd}

\newtheoremstyle{theoremee}
  {6pt}
  {6pt}
  {\itshape}
  {0pt}
  {\bfseries}
  {.}
  { }
  {\thmname{#1}\thmnumber{ #2}\textnormal{\thmnote{ #3}}}
  \theoremstyle{theoremee}

\newcommand{\xrar}[1]{\xrightarrow{#1}}

\newcommand{\riso}{\xrar{\sim}}

 \newenvironment{psmat}
  {\left(\begin{smallmatrix}}
  {\end{smallmatrix}\right)}

 \newenvironment{psmatrix}
  {\left(\begin{smallmatrix}}
  {\end{smallmatrix}\right)}
  
 \newenvironment{pmat}
  {\begin{pmatrix}}
  {\end{pmatrix}}

\newcommand{\textboxnobrace}[2] {\parbox{#1em}{\center{#2}}}

\newcommand{\overtext}[2]{\mathrel{\overset{\makebox{\text{\normalfont\tiny\sffamily#1}}}{#2}}}

\newcommand{\undertext}[2]{\mathrel{\underset{\makebox{\text{\normalfont\tiny\sffamily#1}}}{#2}}}

\newcommand{\wt}{\widetilde}
\newcommand{\wh}{\widehat}

\newcommand{\ov}{\overline}

\newcommand{\toolong}[1]{\resizebox{43em}{!}{$#1$}}

\newcommand{\tr}{\operatorname{tr}}
\newcommand{\rar}{\rightarrow}

\newcommand{\Rar}{\Rightarrow}

\newcommand{\ncom}[1]{\newcommand{#1}}
\ncom{\sbuset}{\subset}

\newcommand{\hrar}{\hookrightarrow}
\newcommand{\thrar}{\twoheadrightarrow}

\DeclareSymbolFont{cyrletters}{OT2}{wncyr}{m}{n}
\DeclareMathSymbol{\Sha}{\mathalpha}{cyrletters}{"58}

\newcommand{\bs}{\backslash}
\newcommand{\et}{\operatorname{\acute{e}t}}

\DMO{\bmr}{\mathbbm{r}}
\DMO{\bmf}{\mathbbm{f}}
\DMO{\bmx}{\mathbbm{x}}

\newcommand{\bA}{\mathbb{A}}

\newcommand{\bC}{\mathbb{C}}

\newcommand{\bG}{\mathbb{G}}
\newcommand{\bH}{\mathbb{H}}

\newcommand{\bN}{\mathbb{N}}

\newcommand{\bQ}{\mathbb{Q}}
\newcommand{\bR}{\mathbb{R}}
\newcommand{\bS}{\mathbb{S}}

\newcommand{\bZ}{\mathbb{Z}}

\newcommand{\sA}{\mathscr{A}}

\newcommand{\sD}{\mathscr{D}}

\newcommand{\sH}{\mathscr{H}}

\newcommand{\sM}{\mathscr{M}}

\newcommand{\sR}{\mathscr{R}}

\newcommand{\cA}{\mathcal{A}}

\newcommand{\cE}{\mathcal{E}}

\newcommand{\cH}{\mathcal{H}}

\newcommand{\cM}{\mathcal{M}}

\newcommand{\cO}{\mathcal{O}}

\newcommand{\cW}{\mathcal{W}}

\newcommand{\fC}{\mathfrak{C}}

\newcommand{\fP}{\mathfrak{P}}

\newcommand{\fX}{\mathfrak{X}}

\newcommand{\fa}{\mathfrak{a}}
\newcommand{\fb}{\mathfrak{b}}

\newcommand{\fg}{\mathfrak{g}}
\newcommand{\fh}{\mathfrak{h}}
\newcommand{\fk}{\mathfrak{k}}

\newcommand{\fn}{\mathfrak{n}}
\newcommand{\fo}{\mathfrak{o}}
\newcommand{\fp}{\mathfrak{p}}
\newcommand{\fq}{\mathfrak{q}}

\newcommand{\fs}{\mathfrak{s}}
\newcommand{\ft}{\mathfrak{t}}

\newcommand{\fsl}{\mathfrak{sl}}

\DMO{\str}{str}
\DMO{\comp}{comp}
\DMO{\rel}{rel}
\DMO{\rat}{rat}
\DMO{\Crys}{Crys}
\DMO{\Pet}{Pet}
\DMO{\Dol}{Dol}
\DMO{\Vect}{Vect}
\DMO{\Ssp}{sp}

\DMO{\Out}{Out}
\DMO{\Inv}{Inv}
\DMO{\MIC}{MIC}
\DMO{\Isoc}{Isoc}
\DMO{\corank}{corank}
\DMO{\Td}{Td}
\DMO{\hol}{hol}
\DMO{\gen}{gen}
\begin{document}
\title{Coherent cohomology of Shimura varieties, motivic cohomology, and archimedean $L$-packets}\author{Gyujin Oh}
\maketitle
\begin{abstract}
We formulate an analogue of the archimedean motivic action conjecture of Prasanna--Venkatesh for \emph{irregular} cohomological automorphic forms on Shimura varieties, which appear on multiple degrees of coherent cohomology of Shimura varieties. Such multiple appearances are due to many infinity types in a single $L$-packet with equal minimal $K$-types. Accordingly, we formulate the conjecture comparing periods of forms of \emph{different automorphic representations}. We provide evidences for the conjecture by showing its compatibility with existing conjectures on periods of automorphic forms. The conjectures suggest the existence of certain operations which move between different infinity types in an $L$-packet. 
\end{abstract}
\tableofcontents
\section{Introduction}
The motivic action conjecture of Venkatesh posits that, roughly speaking, for a Hecke eigensystem $h$, there is a natural action of the motivic cohomology of the adjoint motive associated to $h$ on the $h$-isotypic part of the rational cohomology of locally symmetric spaces. This has many incarnations which shed new light on various parts of the Langlands program. However, the conjectures have been mostly restricted to the case of ``$\delta\ne0$,'' namely when the reductive group in concern $G$ has no compact Cartan subgroup. This in particular excludes the case when the locally symmetric space is a \emph{Shimura variety}. As a locally symmetric space that is not a Shimura variety so far has not been related to algebraic geometry (although see \cite{GoreskyTai}), the motivic action conjectures seemed to be extremely difficult to approach.

On the other hand, there have been expectations that a similar conjecture would exist for automorphic forms over Shimura varieties with \emph{irregular weight}. The easiest instance is the case of weight one modular forms; they appear in both $H^{0}$ and $H^{1}$ of the modular curve of the same line bundle. The main purpose of the paper is to give a formulation of such conjecture for general Shimura varieties, and provide somewhat intricate evidences using well-known results and conjectures regarding periods of automorphic forms. 
The following is a generalization of the archimedean motivic action conjecture to general Shimura varieties.
\begin{customconj}{1}[(Archimedean motivic action conjecture for Shimura varieties)]\label{MainIntro}Let $\lambda$ be a nondegenerate singular analytically integral character, and let $\Pi$ satisfy Assumption \ref{Assumption}, with $\Pi_{\infty}\in\fP_{\lambda}$. Let $\cM=H_{M}^{1}((\Ad\Pi)_{\cO_{E}},\ov{\bQ}(1))$ and $\cH^{i}=H^{i}(X_{G}(\Gamma)_{\ov{\bQ}},[V])[\Pi_{f}]$, where both are regarded as $\ov{\bQ}$-vector spaces equipped with Hermitian bilinear forms, induced from a fixed admissible bilinear form on $\fg_{\bC}$ (see \S\ref{MetricSection}). Then, there is an isometry between graded $\ov{\bQ}$-vector spaces equipped with Hermitian metrics,
\[\wedge^{*}\cM^{*}\otimes\cH^{i_{\min}}\cong \bigoplus_{i=i_{\min}}^{i_{\max}}\cH^{i},\]where $V$ is the automorphic vector bundle coming from Levi (see Notation) such that $\Pi_{f}$ appears in multiple degrees of its cohomology, $i_{\min}=\min\lbrace i\mid H^{i}(X_{G}(\Gamma),[V])[\Pi_{f}]\ne0\rbrace$, and $i_{\max}$ is defined analogously.
\end{customconj}
This is the analogue of \cite[Conjecture 1.2.1]{PV}, although new subtleties arise in the irregular weight case as we will see. Under some standard conjectures on periods, we show that the conjecture indeed holds in a few low-dimensional cases.
\begin{customthm}{2}\label{MainEvidenceIntro}Under certain mild conditions and several standard conjectures on periods (see below), for $G=\Sp_{4}$ and $\SU(2,1)$, Conjecture \ref{MainIntro} is true. In other words, the action of adjoint motivic cohomology group on the coherent cohomology groups of Shimura variety respects $\ov{\bQ}^{\times}$-structure.
\end{customthm}
The ``mild conditions'' are that, informally speaking, the finite part is globally generic, and that there are newforms; see Assumption \ref{Assumption}. These conditions exist to have a clean statement of the conjectures. More importantly, the ``standard conjectures on periods'' are summarized in Assumption \ref{PeriodAssumption}, which includes the Lapid--Mao conjecture and the Beilinson conjectures. The proof of Theorem \ref{MainEvidence} uses the same idea of \cite{PV}, but requires more machinery, as the whole setup is about comparing periods of \emph{different automorphic representations}. 

The cases where Theorem \ref{MainEvidenceIntro} is proved are when the Hecke eigensystem appears in two consecutive degrees of cohomology. In this case, or more generally, when we restrict the statement of Conjecture \ref{MainIntro} into the relation between the top and bottom degrees, the conjecture can be made into a statement that does not refer to motivic cohomology. 

More precisely, let $\Pi=\Pi_{f}\otimes\Pi_{\infty}$ be a cuspidal automorphic representation satisfying Assumption \ref{Assumption}, and let $\lambda$ be the infinitesimal character of $\Pi_{\infty}$, which is singular and nondegenerate (see Notation). Let $V$ be the automorphic vector bundle coming from the Levi (see Notation), such that $H^{i}(\fp,K;V\otimes\Pi_{\infty})\ne0$ for some $i$, and $\Pi_{\min},\Pi_{\max}$ be the members of the archimedean $L$-packet of $\Pi_{\infty}$ such that the degree that $\Pi_{\min}$ ($\Pi_{\max}$, respectively) has nontrivial $(\fp,K)$-cohomology with coefficient in $V$ is the minimum (maximum, respectively) in the $L$-packet. Let $i_{\min}$ ($i_{\max}$, respectively) be the degree, and let $f_{\min}\in\Pi_{f}^{\new}\otimes\Pi_{\min}^{\new}$ and $f_{\max}\in\Pi_{f}^{\new}\otimes\Pi_{\max}^{\new}$ (see Definition \ref{NewDef}) such that $[f_{\min}],[f_{\max}]\in H^{*}(X_{G}(\Gamma),V)$ (the harmonic Dolbeault forms corresponding to the automorphic forms; see Definition \ref{CohclassDef} for a precise definition) are defined over $\ov{\bQ}$. Then, under the Beilinson conjectures, the information on the top and bottom degrees in Conjecture \ref{MainIntro} is precisely that\[\frac{\langle f_{\min},f_{\min}\rangle}{\langle f_{\max},f_{\max}\rangle}\sim_{\ov{\bQ}^{\times}\cap\bR}\left|\pi^{i_{\min}-i_{\max}}\frac{L_{\infty}(1,\Pi,\Ad)}{L_{\infty}(0,\Pi,\Ad)}\cdot\frac{L(1,\Pi,\Ad)}{\vol F^{1}H_{\dR}(\Ad\Pi)}\right|^{2},\]where $\langle\cdot,\cdot\rangle$ is the Petersson inner product, and the volume is computed with respect to the metric induced by any weak polarization (see \cite[\S2.2.3]{PV})\footnote{The volume is independent of the choice of weak polarization, see \cite[Lemma 2.2.2]{PV}.}. In particular, this statement is \emph{equivalent to Conjecture \ref{MainIntro} if there are only two degrees that $\Pi_{f}$ appear in the cohomology}, assuming Beilinson conjectures. 

At first sight, the information on top and bottom degrees might seem to less interesting, as one might guess there is an extra duality that one can compare the top and bottom degrees. Indeed, in \cite{PV}, the top and bottom degrees were complementary, so that they are related via duality. On the other hand, in the irregular setting, there are numerous cases where the top and bottom degrees are not complementary. Rather, they are determined by the position of the infinitesimal character inside the Weyl chambers. Indeed, we work with two examples, $G=\Sp_{4}$ and $\SU(2,1)$, and both cases we work with the choice of $\lambda$ where the minimum and maximum degree of appearances are $H^{0}$ and $H^{1}$, respectively.
\subsection{Archimedean $L$-packet and the matter of choosing automorphic realizations}We now explain the new subtleties of the irregular weight case. 
In view of automorphic cohomology, the phenomenon of a weight one modular form appearing in $H^{0}$ and $H^{1}$ actually involves \emph{two different representations}. If we denote $\omega$ to be the weight one line bundle or the corresponding representation of $\SO(2)$, then for a weight one modular newform $f$,
\[H^{0}(X,\omega)[f]=H^{0}(\fp,\SO(2);D_{0}^{+}\otimes\omega),\]\[H^{1}(X,\omega)[f]=H^{1}(\fp,\SO(2);D_{0}^{-}\otimes\omega),\]where $\fp=\fp_{-}\oplus\fs\fo(2)$, $\fp_{-}$ is the anti-holomorphic tangent space and $D_{0}^{+}$ and $D_{0}^{-}$ are the holomorphic and the antiholomorphic discrete series, respectively. Indeed, it is $\ov{f}d\ov{z}$ that appears in $H^{1}(X,\omega)$ as a class in Dolbeult cohomology, and $\ov{f}$ is an \emph{antiholomorphic} modular form of weight $-1$. 

In general, if the same Hecke eigensystem appears in multiple degrees of cohomology of the same automorphic vector bundle of a Shimura variety, then each such instance is actually represented by the so-called $(\fp,K)$-cohomology of different automorphic representations. More precisely, the finite part ($G(\bA_{f})$-representation) remains the same, while the infinity type varies inside an archimedean $L$-packet. Such phenomenon happens precisely when the infinitesimal character of the infinity type lies on the walls of Weyl chambers. 

The action of $\wedge^{*}\fa_{G}^{*}$ in \cite{PV} is most naturally thought as the self-$\Ext$-algebra of a representation. On the other hand, in our case, the representations corresponding to the target and the source of the action is different, so the action cannot be thought as an $\Ext$-\emph{algebra}, but merely as an $\Ext$-group. Furthermore, the translation of the action into the automorphic cohomology context also depends on the choice of \emph{automorphic realization maps $\Pi_{f}\otimes\Pi_{\infty}\hrar\cA(G)$} for each $\Pi_{\infty}$; namely, each realization map can be always scaled by a scalar, but we want to compare them as a whole. This problem did not arise in \emph{op. cit.}, as a choice of a single automorphic realization map would rigidify the situation. In turn, we had no choice but to formulate a slightly weaker Conjecture \ref{MainIntro} that asserts the motivic action conjecture on the level of metric spaces. 

\subsection{Comparison of periods of different automorphic representations}
We briefly explain the strategy of the proof of Theorem \ref{MainEvidenceIntro}. 
As does in \cite{PV}, we most notably assume \emph{Beilinson's conjectures (for Chow motives)}. The main new feature in this paper is that, because we need to compare periods of two different representations, we need two different conjectures on periods, one for each representation. In both evidences, we will compare periods of a \emph{holomorphic automorphic form} and a \emph{generic automorphic form}. For the holomorphic form, we will need the \emph{refined Gan--Gross--Prasad conjecture} (referred as Ichino--Ikeda conjecture in \emph{op. cit.}), although in some instances this requirement can be avoided by using the doubling construction of standard $L$-functions. For the generic form, we will need the \emph{Lapid--Mao conjecture}, which relates the value of a Whittaker function to a certain $L$-value.

One also needs a way to detect rationality of classes for both types of forms. For the holomorphic forms, we can use Fourier expansion, but for those appearing in higher coherent cohomology, we need a new machinery. We will develop the so-called \emph{cohomological period integrals} for higher coherent cohomology, which realizes integral representations of $L$-functions as cup product pairings in coherent cohomology. The cohomological interpretation of such integral representations (or, at least, their appearance in the literature) is relatively new (\cite{LPSZ}, \cite{Oh}).\subsection{Generalized complex conjugations}
The archimedean motivic action conjecture as stated in \cite{PV} gives a recipe of rational cohomology classes, whereas our Conjecture \ref{MainIntro} is a statement on metrics. To formulate a similar conjecture in the setting of coherent cohomology of Shimura varieties, we need a way to rigidify between different automorphic representations in an archimedean $L$-packet. Indeed, in retrospect, even in the easy case of modular forms, one needs complex conjugation to go between holomorphic and antiholomorphic limits of discrete series. Unfortunately, beyond the case of modular forms, there is no known general operation that can move between different infinity types. We will tentatively name such an operation a \ \emph{generalized complex conjugation}, which should send an automorphic form of a certain infinity type to an automorphic form of another infinity type in the same $L$-packet. It seems inevitable to come up with such an operation to formulate the full conjecture on rational cohomology classes.

The generalized complex conjugations should be naturally understood in the context of a ``derived'' local-global compatiblity in some sense, and their existence is also suggested by the existence of similar operations in the analogous settings over the $p$-adic fields (Kottwitz's conjectures, e.g. \cite{FarguesMantovan}) and over the function fields (excursion operators, e.g. \cite{Lafforgue}). Following the suggestions of Joseph Wolf, we will investigate the nature of generalized complex conjugations using the theory of Penrose transforms.

On the other hand, if the associated Hermitian symmetric space is a product of copies of the upper half planes, one can come up with an operation that changes one infinity type to another by taking complex conjugation at certain variables. This is a \emph{partial complex conjugation}, studied in \cite{HarrisHilb}. Using partial complex conjugations, we can formulate a conjecture on rationality of cohomology classes in the case of Hilbert modular forms of partial weight one. There exists a prior work of \cite{Horawa} on the motivic action conjecture for Hilbert modular forms of parallel weight one, which similarly uses partial complex conjugations. We compare our conjecture in the Hilbert modular form case with the conjecture of \cite{Horawa}, and explain the evidences given in \emph{op. cit.} are also consistent with our conjecture.

\subsection{Summary}
In \S\ref{Sectiontwo}, we take an efficient route to the statement of the Archimedean motivic action conjecture (Conjecture \ref{Metric}) and its more accessible variant, the Period conjecture (Conjecture \ref{Period}). The objective of the section is to set the conjecture in a context. We in particular defer the abstract discussion of how to derive the conjectures, parallel to those of \cite[\S2-\S5]{PV}, to later sections, as it requires more advanced theory on representation theory of real groups.

In \S\ref{Sp4EvidenceSection} and \S\ref{SU21EvidenceSection}, we provide our main evidence for the Period conjecture (Conjecture \ref{Period}). Similarly to \cite[\S7]{PV}, we prove that, for $G=\Sp_{4}$ and $\SU(2,1)$, the Period conjecture is compatible with several well-accepted conjectures on periods of automorphic forms, such as the Beilinson's conjectures, the Lapid--Mao conjecture and the refined Gan--Gross--Prasad conjectures. To use these, we review how certain period integrals can be interpreted as cup product pairings of (higher) coherent cohomology classes on Shimura varieties.

In \S\ref{MotivicRationalSect}, we discuss the issue on formulating a motivic action conjecture on rationality of coherent cohomology classes. Most notably, we suggest the notion of \emph{generalized complex conjugations}, which move between different members of a single archimedean $L$-packet. In \S\ref{HorawaSect}, we formulate a precise conjecture in the case of Hilbert modular forms using partial complex conjugations, and compare our conjecture with the conjecture of \cite{Horawa}. In \S\ref{GCCSection}, we spell out conditions that the generalized complex conjugations should satisfy, and formulate the full conjecture assuming their existence. In Appendix \ref{BCSection}, we review the formulation of Beilinson's conjecture for motives over a general number field, as many references state the conjecture for only $\bQ$-motives. Finally, we develop a representation theory background in Appendix \ref{ExtSection}, parallel to \cite[\S2-\S4]{PV}. Although Appendix \ref{ExtSection} is independent of the development of the rest of the paper, the section is suggestive of a correct foundation in which the motivic action conjectures need to be developed.
\subsection{Problems and questions}There are several interesting questions that arise in this work.
\begin{enumerate}
\item Place the generalized complex conjugations in the context of some form of ``derived'' local-global compatibility, motivated by the strong form of Arthur conjectures as realized in the function field case via excursion operators as in \cite{Lafforgue}. A correct formulation should be in accordance with the existing statements of derived local-global compatibility as in \cite{Feng} and \cite{Zhu}. 
\item The compatibility between the Period Conjecture, Conjecture \ref{Period}, and the existing period conjectures relies on the yet-to-be-calculated archimedean zeta integrals. These can be conducted using explicit integral formulae of (generalized) Whittaker functions, e.g. \cite{KO}, \cite{OdaSp}.
\item As we deal with motives over a more general number field and Shimura varieties over a number field other than $\bQ$, in every aspect of our discussion, the choice of a complex embedding is always implicit. In particular, there must be a relation between the conjectures in this paper for the \emph{conjugates of Shimura varieties} (e.g. \cite{Var}).
\item It seems extremely hard to detect rationality of coherent cohomology classes if the Hecke operators can only cut a space that is of dimension larger than one. For example, if $X$ is a Hilbert modular surface, and if $\omega$ is the parallel weight one line bundle, then it seems extremely difficult to determine whether a class in $H^{1}(X,\omega)$ is defined over $\ov{\bQ}$, or even to produce a class in it. 
\item It is expected that the motivic action conjecture will involve $L$-packets even in the case of $\delta>0$. It may be possible to formulate the conjecture for the same eigensystem appearing in cohomology with \emph{different coefficients}.
\end{enumerate}

\subsection{Notation}\label{Notation}
Let $G$ be a connected reductive algebraic group over $\bQ$. For simplicity, let us assume that $G$ is quasisplit, $G(\bR)$ is connected, and that the center of $G$ does not have a nontrivial $\bR$-split torus. Also, we assume that there exists a \emph{twisting element} in the sense of \cite[Definition 5.2.1]{BuzzardGee}\footnote{This is to avoid the subtlety of difference between $C$-algebraicity and $L$-algebraicity.}. Let $\fg_{\bQ}$ be the $\bQ$-Lie algebra of $G$, and let $\fg_{\bR},\fg_{\bC}$ be its base change to $\bR$ and $\bC$, respectively. We occasionally drop the subscript for $\bC$. Let $W_{G}$ be the Weyl group of $G$. We endow an invariant, $\theta$-invariant, $\bR$-invariant bilinear form $B$ on $\fg_{\bR}$, such that $B(X,\theta(X))$ is negative definite, where $\theta$ is the Cartan involution. We will use this to talk about inner product on weight space, Riemannian metric on the Hermitian symmetric domain, etc. In specific examples, we may and will choose $B$ to induce a preferred Riemannian metric on the Hermitian symmetric space (for example, one may want the Riemannian metric to be $dxdy$ on the upper half plane $\bH=\lbrace x+iy\mid y>0\rbrace$). We also use $B$ to any other bilinear form induced from $B$. 

Suppose further that $G$ gives rise to a \emph{Shimura variety}, which means that there is a symmetric space $X$ for $G(\bR)$ which can be endowed with a structure of Hermitian symmetric domain (which we will fix). Fix a point $h\in X$, which gives rise to a \emph{Hodge cocharacter} $h:\bS=\Res_{\bC/\bR}\bG_{m,\bC}\rar G_{\bR}$ which in turn induces a real Hodge structure of weight $0$ on $\fg_{\bR}$, $\fg=\fg^{-1,1}\oplus\fg^{0,0}\oplus\fg^{1,-1}$. Given an open compact subgroup $\Gamma\subset G(\bA_{f})$, there exists a quasi-projective variety $Y_{G}(\Gamma)$, a Shimura variety, defined over a number field $E$, whose complex points\footnote{Note that the choice of a point in a Hermitian symmetric domain gives the reflex field as a subfield of $\bC$, so there is a \emph{preferred complex embedding}; e.g. \cite[Notation 4.6]{Var}. In particular, one can expect that the statement of the conjecture depends a priori on the choice of a Hermitian symmetric domain. It could be interesting to check if our conjecture is consistent with conjugation of Shimura varieties.\label{reflexfoot}} have an analytification isomorphic to the double quotient $G(\bQ)\bs (X\times G(\bA_{f})/\Gamma)$. 

Let $K\subset G(\bR)$ be the stabilizer of $h$. Let $T$ be the Cartan subgroup of $K$. Then, $X\cong G(\bR)/K$ and $\fg^{0,0}=\fk:=\Lie(K)_{\bC}$. We denote $\fp_{+}=\fg^{-1,1}$, $\fp_{-}=\fg^{1,-1}$ and $\fp=\fk\oplus\fp_{-}$. Then, $\fp$ is a parabolic subalgebra of $\fg$, giving rise to a parabolic subgroup $P\subset G_{\bC}$ with $\Lie P=\fp$. We also fix once and for all a positive system of roots for $\fk$. The holomorphic tangent space of $X$ at $h$ is identified with $\fp_{+}$, so there is a $G(\bR)$-equivariant embedding of complex manifolds $X\rar \check{D}:=G(\bC)/P(\bC)$, sending $h\mapsto P(\bC)$. In this regard, $K=G(\bR)\cap P(\bC)$, and $K(\bC)$ is the Levi subgroup of $P(\bC)$. Also, any finite-dimensional holomorphic representation $V$ of $P(\bC)$ gives rise to an algebraic vector bundle over $Y_{G}(\Gamma)$, an \emph{automorphic vector bundle}, denoted $[V]$, which is an algebraization of the pullback of the holomorphic vector bundle on $X$ which in turn is the restriction of the vector bundle $G(\bC)\times^{P(\bC)}V\rar \check{D}$ on $\check{D}$. If $V$ factors through $K(\bC)$, namely if it is induced from a representation of $K(\bC)$, we will call $[V]$ an automorphic vector bundle \emph{coming from the Levi}. If not, we will call $V$ \emph{nearly}, following \cite{LPSZ}.

To save space, we may abbreviate some words with repeated appearances: discrete series into DS, limit of discrete series into LDS, and nondegenerate limit of discrete series (see the paragraph before Theorem \ref{PK} for its definition) into NLDS. All real group representations are thought as $(\fg,K)$-modules.

For automorphic representations, their $L$-functions are normalized so that $\frac{1}{2}$ is the center of symmetry. For pure motives, theire $L$-functions are normalized so that, if $w$ is its weight, $\frac{w+1}{2}$ is the center of symmetry. For both kinds of $L$-functions, $w$ is called the \emph{motivic weight} of the $L$-function. For an automorphic representation $\Pi=\Pi_{f}\otimes\Pi_{\infty}$ of $G(\bA)$, the field of rationality $F_{\Pi}$ is the fixed field of the isomorphism class of $\Pi_{f}$ as a $G(\bA_{f})$-representation (\cite[\S3.1]{Clozel}). An inner product on the space of automorphic forms, denoted $\cA(G)$, can be given as either the $L^{2}$-norm on $G(\bQ)/G(\bA)$ with respect to the Tamagawa measure or the measure coming from the Riemannian metric of the symmetric space, as the norms are all scaled by the same scalar factor. The second norm is the same as the usual \emph{Petersson norm}, which we will denote as $\langle ,\rangle_{P}$.

For the integral representations, we fix a nontrivial additive character $\psi=\prod_{p}\psi_{p}$ of $\bA/\bQ$. For a cuspidal automorphic form $\varphi=\prod_{p}\varphi_{p}$ of cuspidal automorphic representation $\pi=\bigotimes_{p}'\pi_{v}$ of $G(\bA)$, its Whittaker transform $W_{\varphi}(g)$ is defined as $W_{\varphi}(g)=\int_{N_{G}(\bQ)\bs N_{G}(\bA)}\varphi(ng)\psi(n^{-1})dn$, for some choice of $N_{G}$ that needs to be specified when talking about Whittaker model. Locally, a generic $G(\bQ_{p})$-representation $\pi_{p}$ is isomorphic to the space of functions $\cW_{p}:=\lbrace W_{p}:G(\bQ_{p})\rar\bC\mid W_{p}(ng)=\psi_{p}(n)W_{p}(g)\text{ for }n\in N_{G}(\bQ_{p})\rbrace$; one can choose isomorphisms (a local \emph{Whittaker model}) $\pi_{p}\riso\cW_{p}$, $f_{p}\mapsto W_{f_{p}}$ that are compatible with the global Whittaker model; namely, $W_{\varphi}(g)=\prod_{p}W_{\varphi_{p}}(g_{p})$. 

We use notational convention for motives as in \cite[\S2]{PV}, except that we will deal with motives over a more general number field. In that case, we put the complex embedding in the subscript, such as $M_{\sigma}$, $\comp_{B,\dR,\sigma}$, etc.
\section{Archimedean $L$-packets and the motivic action conjecture}\label{Sectiontwo}
In this section, we take the shortest path to the statement of the Archimedean motivic action conjecture for Shimura varieties, Conjecture \ref{Metric}. More abstract justification of the formulation of the Conjecture, including a parallelism between \cite{PV} and our conjecture, is discussed in Appendix \ref{ExtSection} and \S\ref{GCCSection}.
\subsection{$(\fp,K)$-cohomology and automorphic forms}
Firstly, we quickly review how coherent cohomology of Shimura varieties is related to automorphic forms via the theory of $(\fp,K)$-cohomology. Recall that, as the singular cohomology of locally symmetric spaces can be calculated in terms of $(\fg,K)$-cohomology, the coherent cohomology of Shimura varieties can be calculated in terms of the so-called $(\fp,K)$-cohomology. By reinterpreting what the Dolbeault cohomology calculates in the setting of Shimura varieties, one gets the following
\begin{prop2}[(See {\cite[(2.12)]{Jun}})]We have
\[H^{i}(Y_{G}(\Gamma),[V])\cong H^{i}(\fp,K;C^{\infty}(G(\bQ)\bs G(\bA)/\Gamma)^{K\mathrm{-finite}}\otimes V),\]where the left hand side is analytic cohomology, and $V$ is understood as a $(\fp,K)$-module with trivial $\fp$-action. 
\end{prop2}
Furthermore, there is an analogue of Franke's theorem for coherent cohomology of Shimura varieties.
\begin{thm2}[(Su, {\cite[Theorem 6.7]{Jun}})]\label{Jun}
For any sufficiently refined polyhedral cone decomposition $\Sigma$, there is a natural Hecke-equivariant isomorphism
\[H^{i}(X_{G}^{\Sigma}(\Gamma),[V]^{\can})\cong H^{i}(\fp,K;\cA(G)^{\Gamma}\otimes V),\]where $X_{G}^{\Sigma}(\Gamma)$ is the corresponding toroidal compactification, $[V]^{\can}$ is the canonical extension of $[V]$, and $\cA(G)$ is the space of automorphic forms, namely the space of right $K$-finite, $Z(\fg)$-finite smooth functions on $G(\bQ)\bs G(\bA)$ of moderate growth.
\end{thm2}\begin{rmk2}
We would be only interested in a part of coherent cohomology localized at a cuspidal Hecke eigensystem, so it is unlikely that the full power of Su's theorem is required.
\end{rmk2}The calculation of coherent cohomology is therefore about $(\fp,K)$-cohomology of automorphic representations. As far as the coherent cohomology is concerned, the choice of $\Sigma$ is ineffective, so we may occasionally drop the superscript $\Sigma$ if there is no confusion.

We will be interested in the situation where a Hecke eigensystem appears in multiple degrees of coherent cohomology. Namely, for an admissible $G(\bA_{f})$-representation $\Pi_{f}$ and an automorphic vector bundle $\cE$ of $Y_{G}(\Gamma)$, we are interested in 
\[H^{*}(X_{G}(\Gamma),\cE^{\can})[\Pi_{f}].\]By Theorem \ref{Jun}, we have a canonical isomorphism
\[H^{*}(X_{G}(\Gamma),\cE^{\can})\cong\left(\bigoplus_{\Pi_{f}\otimes\Pi_{\infty}\subset\cA(G)}H^{*}(\fp,K;\Pi_{\infty}\otimes E)\right)\otimes_{\bC}\Pi_{f}^{\Gamma},\]where the sum runs over all automorphic representations with the finite part being $\Pi_{f}$, and $E$ is the algebraic $P(\bC)$-representation such that $[E]=\cE$. Thus, it is possible that \emph{several different $\Pi_{\infty}$'s can appear in the decomposition}, if $H^{*}(\fp,K;\Pi_{\infty}\otimes E)\ne0$ for several different $\Pi_{\infty}$'s. Indeed, this can be the case if, for example, some of $\Pi_{\infty}$ is a nondegenerate limit of discrete series (NLDS) and not a discrete series (DS); recall that a \emph{nondegenerate limit of discrete series} is a limit of discrete series whose infinitesimal character is not orthogonal to any compact root. Unlike the case of ``$\delta>0$'' as in \cite{PV}, the appearance of single Hecke eigensystem in multiple cohomological degrees in our setting necessarily implies that, by the following Theorem, there are \emph{many different archimedean representations involved}:
\begin{defthm2}[(See {\cite{VZ}}, {\cite{Schmid}})]\label{PK}Let $\Pi_{\infty}$ be the $(\fg,K)$-module associated to a DS or a NLDS representation of $G(\bR)$. Then, there is a unique $0\le i\le \dim X$ and a finite-dimensional irreducible $K$-representation $V$ such that $H^{i}(\fp,K;\Pi_{\infty}\otimes V)\ne0$. Furthermore, $\dim_{\bC}H^{i}(\fp,K;\Pi_{\infty}\otimes V)=1$. We will denote $i_{\Pi_{\infty}}$ and $V_{\Pi_{\infty}}$ for the $i$ and $V$ corresponding to $\Pi_{\infty}$.
\end{defthm2}Indeed, the above Theorem says that a single archimedean representation can only contribute to a single degree. We will see that the \emph{raison d'\^etre} of appearance of a Hecke eigensystem in multiple degrees is that the infinity type $\Pi_{\infty}$ can change in an archimedean $L$-packet without changing the finite part. We will see in detail in Appendix \ref{ExtSection} how an archimedean $L$-packet (rather than a single $G(\bR)$-representation) appears in the context of motivic action. For now, we move on to the formulation of the ``metric'' conjecture which does not involve abstract real group representation theory nor Lie algebra cohomology. From now on, for the sake of simplicity, we assume the following
\begin{ass2}\label{Assumption}
Let $\Pi=\Pi_{f}\otimes\Pi_{\infty}$, with $\Pi_{f}=\prod_{p<\infty}\Pi_{p}$, be a cuspidal automorphic automorphic representation with $\Pi_{\infty}$ an NLDS (see Notation)\footnote{This would ensure that $F_{\Pi}$, the field of definition, is a number field, by our assumption in Notation that there exists a twisting element in the sense of \cite[Definition 5.2.1]{BuzzardGee}. Note that the twisting element indeed exists for $G=\Sp_{4}$ (as it is split and has simply-connected derived subgroup) and $\SU(2,1)$ (as the half-sum of positive roots is integral).}. We hereafter assume the following:
\[\label{GEN}\text{$\Pi_{f}$ is globally generic.}\tag{GEN}\]
\[\label{NEW}\text{For each $p<\infty$, there exists a compact open subgroup $\Gamma_{p}\le G(\bZ_{p})$ such that $\dim_{\bC}\Pi_{p}^{\Gamma_{p}}=1$.}\tag{NEW}\]
\[\label{Sp4A}\textboxnobrace{40}{If $G=\Sp_{4}$, a holomorphic Siegel modular newform $f_{\Pi}$ in $\Pi_{f}\otimes\Pi_{\infty}^{\hol}$, where $\Pi_{\infty}^{\hol}$ is a holomorphic (limit of) discrete series, has a nontrivial special Bessel period $B(f_{\Pi},F)\ne0$ for some imaginary quadratic field $F$. Here, $B(f_{\Pi},F)=\sum_{A=\begin{psmat}a&b/2\\b/2&c\end{psmat},4ac-b^{2}=D_{F},A\sim M^{T}AM \text{ for }M\in\SL_{2}(\bZ)}a_{A}$, where $f_{\Pi}(Z)=\sum_{S}a_{S}e^{2\pi i \tr(SZ)}$ is the $q$-expansion of $f_{\Pi}$, and $-D_{F}$ is the discriminant of $F$.}\tag{$\Sp_{4}$}\]
\[\label{SU21A}\textboxnobrace{40}{If $G=\SU(2,1)$, defined using an imaginary quadratic field $K$, (1) there exists a split place $v$ such that $\Pi_{v}$ is supercuspidal, and (2) if $p$ is not split in $K$ where $\Pi_{p}$ is also not unramified, either $\Pi_{p}$ is supercuspidal or the stabilizer of an anisotropic vector $\SU(1,1)(\bQ_{p})\subset\SU(2,1)(\bQ_{p})$ is compact.}\tag{$\SU(2,1)$}\]
\end{ass2}
The purpose of Assumption \ref{Assumption} is to isolate the effect of archimedean $L$-packet phenomenon amongst others and to use existing results on the refined Gan--Gross--Prasad conjectures. We believe that, for example, it would be not difficult to formulate the Conjectures without (\ref{NEW}). 
\begin{defn2}\label{NewDef}Under Assumption \ref{Assumption}, let $\Gamma(\Pi_{f})=\Gamma:=\prod_{p}\Gamma_{p}$. Also, when we say a vector $f\in\Pi$ is a \emph{newform}, it means $f=\prod_{p}f_{p}\otimes f_{\infty}$ such that, not only $f_{p}\in\Pi_{p}^{\Gamma_{p}}$, but also $f_{\infty}$ is a highest weight vector of the minimal $K$-type of $\Pi_{\infty}$ (which exists as we will be only concerned about either DS or LDS). We define $\Pi_{f}^{\new}=\prod_{p}\Pi_{p}^{\Gamma_{p}}$, and also $\Pi_{\infty}^{\new}$ to be the one-dimensional $\bC$-vector subspace of $\Pi_{\infty}$ generated by highest weight vectors of the minimal $K$-type. In particular, $f$ being a newform means that $f\in\Pi_{f}^{\new}\otimes\Pi_{\infty}^{\new}$.

In view of Theorem-Definition \ref{PK}, the above ``newform'' appears in 
\[H^{i_{\Pi_{\infty}}}(X_{G}(\Gamma),[V_{\Pi_{\infty}}]^{\can})[\Pi_{f}],\]which we denote by $H^{i_{\Pi_{\infty}}}(X)[\Pi_{f}]$.
\end{defn2}
Note that a choice of a highest weight vector of $V_{\Pi_{\infty}}$ induces a natural isomorphism \[H^{i_{\Pi_{\infty}}}(\fp,K;\Pi_{\infty}\otimes V_{\Pi_{\infty}})\cong\Pi_{\infty}^{\new},\]defined as follows. Given $v\in\Pi_{\infty}^{\new}$, define a $K$-homomorphism \[f:\left(\wedge^{i_{\Pi_{\infty}}}(\fp/\fk)\right)\otimes V_{\Pi_{\infty}}^{*}\rar \Pi_{\infty},\]by sending the highest weight vector (induced from the choice of a highest weight vector of $V_{\Pi_{\infty}}$ and the roots of $\fg$) of the highest $K$-type of the source to $v$ and sending all other $K$-types to zero. This defines a class in $\Hom_{K}(\wedge^{i_{\Pi_{\infty}}}(\fp/\fk),\Pi_{\infty}\otimes V_{\Pi_{\infty}})$ which is closed in the corresponding Chevalley--Eilenberg complex for the $(\fp,K)$-cohomology, thus a class in $H^{i_{\Pi_{\infty}}}(\fp,K;\Pi_{\infty}\otimes V_{\Pi_{\infty}})$. \begin{rmk2}In this paper, each conjecture will consider a fixed finite type $\Pi_{f}$ and a coefficient vector bundle $V$, so in that context we firstly fix a choice of a highest weight vector of $V$ before anything else.\end{rmk2}
\begin{defn2}\label{CohclassDef}
For $f=\prod_{p\le\infty}f_{p}\in\Pi_{f}^{\new}\otimes\Pi_{\infty}^{\new}$, we define $[f]\in H^{i_{\Pi_{\infty}}}(X)[\Pi_{f}]$ to be the class corresponding to\[[f_{\infty}]\otimes\prod_{p<\infty}f_{p}\in H^{i_{\Pi_{\infty}}}(\fp,K;\Pi_{\infty}\otimes V_{\Pi_{\infty}})\otimes\Pi_{f}^{\new}\subset H^{i_{\Pi_{\infty}}}(X)[\Pi_{f}],\]where $[f_{\infty}]\in H^{i_{\Pi_{\infty}}}(\fp,K;\Pi_{\infty}\otimes V_{\Pi_{\infty}})$ corresponds to $f_{\infty}\in\Pi_{\infty}^{\new}$ via the above natural isomorphism $H^{i_{\Pi_{\infty}}}(\fp,K;\Pi_{\infty}\otimes V_{\Pi_{\infty}})\cong \Pi_{\infty}^{\new}$.
\end{defn2}
\subsection{Metrics on cohomology}\label{MetricSection}
To formulate the main Conjecture \ref{Metric}, we need to define metrics on the coherent cohomology of Shimura varieties as well as motivic cohomology. To define metrics on both motivic cohomology and Dolbeault cohomology, we fix an \emph{admissible bilinear form} on $\fg_{\bC}$ in the following sense.
\begin{defn2}[(Admissible bilinear form)]An \emph{admissible bilinear form} $B$ on $\fg_{\bC}$ is an invariant, $\theta$-invariant Hermitian bilinear form on $\fg_{\bC}$ such that the following conditions are satisfied.
\begin{enumerate}
\item It is a natural extension of a $\bR$-valued bilinear form on $\fg_{\bR}$ such that $B(X,\theta(X))$ is negative definite, where $\theta$ is the Cartan involution that fixes $\fk_{\bR}$.
\item It is $\ov{\bQ}$-valued on $\fg_{\ov{\bQ}}$.
\end{enumerate}
\end{defn2}
For example, if $G(\bR)$ is semisimple, the extention of the Killing form as a Hermitian bilinear form on $\fg_{\bC}$ is an admissible bilinear form.

Firstly, the coherent cohomology of toroidal compactifications of Shimura varieties \[H^{*}(X_{G}(\Gamma),[V]^{\can}),\] can be given a Hermitian metric\footnote{This coincides with the metric on the Lie algebra cohomology (see e.g. \cite[\S II.2]{BorelWallach}), and this will be reviewed later in Appendix \ref{ExtSection}.}, induced from a Hermitian metric on the Hermitian symmetric domain and the automorphic vector bundle (which is always possible by the compactness of $K$, see e.g. \cite[\S5]{Gil}), which is in turn induced from our choice of admissible bilinear form on $\fg_{\bC}$. The choice of Hermitian metric induces Hermitian metrics on the entries of Dolbeault complex $\sA^{0,*}([V]^{\can})$. By taking the formal adjoint $\ov{\partial}^{*}$ of $\ov{\partial}$, one can define the Laplacian $\Delta=\ov{\partial}\ov{\partial}^{*}+\ov{\partial}^{*}\ov{\partial}$ on each entry of Dolbeault complex. By Hodge theory, the Dolbeault cohomology $H^{i}([V]^{\can})$ is identified with the space of \emph{harmonic $(0,i)$-forms} $\sH^{i}([V]^{\can})$, which is just the kernel of the Laplacian $\Delta$. The restriction of Hermitian metric on $\sA^{0,i}([V]^{\can})$ to the space of harmonic $(0,i)$-forms gives rise to a Hermitian metric on the Dolbeault cohomology.

On the other hand, the motivic cohomology that would have to appear in the motivic action conjectures is that of the adjoint motive. As in \cite[\S4.2]{PV}, we assume a conjecture on the existence of adjoint motive.\begin{conj2}\label{AdjointMotive}For $\Pi=\Pi_{f}\otimes\Pi_{\infty}$ satisfying Assumption \ref{Assumption}, there exists an adjoint motive $\Ad\Pi$, in the sense of \cite[Definition 4.2.1]{PV}, over the reflex field $E$ with coefficients in $F_{\Pi}$.\end{conj2}To endow a Hermitian metric on the adjoint motivic cohomology $H^{1}_{M}((\Ad\Pi)_{\cO_{E}},\ov{\bQ}(1))$, consider the Beilinson regulator for $\Ad\Pi$. Recall that the Beilinson regulator is a map from motivic cohomology to Deligne cohomology,
\[H^{1}_{M}((\Ad\Pi)_{\cO_{E}},\bQ(1))\rar H_{\sD}^{1}((\Ad\Pi)_{\bR},\bR(1)),\]where the Deligne cohomology group of a motive $M$ over a number field $k$ is defined as (see \cite[(6.1.22)]{Ramakrishnan})
\[H_{\sD}^{i}(M_{\bR},A)=\prod_{w\text{ complex places of $k$}}H_{\sD}^{i}(M\times_{k,w}\bC,A)\times\prod_{w\text{ real places of $k$}}H_{\sD}^{i}(M\times_{k,w}\bR,A).\]Under the Beilinson's conjectures, the Beilinson regulator gives rise to an isomorphism \[H^{1}_{M}((\Ad\Pi)_{\cO_{E}},\bQ(1))\otimes_{\bQ}\bR\riso H_{\sD}^{1}((\Ad\Pi)_{\bR},\bR(1)).\]
\begin{rmk2}
Note that the Betti realization, on which the Beilinson regulator depends, \emph{depends on the complex embedding} $E\hrar\bC$. In our discussion, we use the preferred embedding that came with the datum of reflex field\footref{reflexfoot}. In the following discussions, we always use this complex embedding.
\end{rmk2}
From \S\ref{DeligneExtSection}, we know that the target of the Beilinson regulator is identified with $\wh{\fg}^{\varphi(W_{\bC/\bR})}$, where $\varphi:W_{\bC/\bR}\rar {}^{L}G$ is the corresponding Lanvlands parameter, and that this is identified as a Lie subalgebra of $\wh{\ft}$. We define a Hermitian bilinear form on $\wh{\ft}$ as the dual Hermitian bilinear form of the one we chose for $\ft$, and this restricts to a Hermitian bilinear form on the Deligne cohomology. 
\begin{defn2}\label{HermitianQbar}
Let $X$ be a finite-dimensional $\ov{\bQ}$-vector space, together with an embedding $\iota:\ov{\bQ}\hrar\bC$. A \emph{Hermitian bilinear form} on $X$ is a Hermitian metric on $X\otimes_{\ov{\bQ},\iota}\bC$.
\end{defn2}By \cite[Lemma 2.2.2]{PV}, the volume of $H_{M}^{1}((\Ad\Pi)_{\cO_{E}},\ov{\bQ}(1))$ is in fact independent of choice of the admissible Hermitian bilinear form.
\subsection{Archimedean motivic action conjecture for Shimura varieties}
Now, we are able to state the Archimedean motivic action conjecture for Shimura varieties as follows.\begin{conj2}[(Archimedean motivic action conjecture for Shimura varieties, metric version)]\label{Metric}
Let $\Pi=\Pi_{f}\otimes\Pi_{\infty}$ satisfy Assumption \ref{Assumption}, with $\Pi_{\infty}$ being a NLDS but not being a DS. Let $\cM=H_{M}^{1}((\Ad\Pi)_{\cO_{E}},\ov{\bQ}(1))$ and $\cH^{i}=H^{i}(X)[\Pi_{f}]$, where both are regarded as $\ov{\bQ}$-vector spaces equipped with a Hermitian bilinear form (see Definition \ref{HermitianQbar}), induced from a fixed admissible Hermitian bilinear form on $\fg_{\bC}$. Then, there is an isomorphism of graded $\ov{\bQ}$-vector spaces equipped with Hermitian metrics,
\[\wedge^{*}\cM^{*}\otimes\cH^{i_{\min}}\cong \bigoplus_{i=i_{\min}}^{i_{\max}}\cH^{i},\]where $i_{\min}$ and $i_{\max}$ are the bottom and top degrees, respectively, of appearance of $\Pi_{f}$ in the cohomology $H^{*}(X)[\Pi_{f}]$.
\end{conj2}\begin{rmk2}
It seems that, to descend the coefficient field from $\ov{\bQ}$ to a number field, one may have to take a field larger than $EF_{\Pi}$ (the compositum of the reflex field and the field of definition of $\Pi_{f}$) even in the case of $\SL_{2}(\bQ)$; see \cite[Corollary 4.6]{Horawa}. \end{rmk2}
The relationship between the above Conjecture and the philosophy of motivic action conjectures will be fully discussed in Appendix \ref{ExtSection}.

A special subset of the Main conjecture (Conjecture \ref{Metric}) concerning the norms of \emph{top and bottom degrees} can be formulated without reference to motivic cohomology, assuming Beilinson's conjectures for Chow motives \cite[Conjecture 2.1.1]{PV}.
\begin{conj2}[(Comparison of top and bottom in Conjecture \ref{Metric})]\label{Period}
Let $\Pi$ be as in Conjceture \ref{Metric}. Let $V$ be the automorphic vector bundle coming from the Levi, such that $H^{i}(\fp,K;V\otimes\Pi_{\infty})\ne0$ for some $i$, and $\Pi_{\min},\Pi_{\max}$ be the members of the archimedean $L$-packet of $\Pi_{\infty}$ such that the degree that $\Pi_{\min}$ ($\Pi_{\max}$, respectively) has nontrivial $(\fp,K)$-cohomology with coefficient in $V$ is the minimum (maximum, respectively) cohomological degree, denoted $i_{\min}$ ($i_{\max}$, respectively) in the $L$-packet. Let $f_{\min}\in\Pi_{f}^{\new}\otimes\Pi_{\min}^{\new}$ and $f_{\max}\in\Pi_{f}^{\new}\otimes\Pi_{\max}^{\new}$ such that $[f_{\min}],[f_{\max}]\in H^{*}(X)[\Pi_{f}]$ are defined over $\ov{\bQ}$. Then, \[\frac{\langle f_{\min},f_{\min}\rangle_{P}}{\langle f_{\max},f_{\max}\rangle_{P}}\sim_{\ov{\bQ}^{\times}\cap\bR}\left|\pi^{i_{\min}-i_{\max}}\frac{L_{\infty}(1,\Pi,\Ad)}{L_{\infty}(0,\Pi,\Ad)}\cdot\frac{L(1,\Pi,\Ad)}{\vol F^{1}H_{\dR}(\Ad\Pi)}\right|^{2},\]where the volume is computed with respect to the metric induced by any weak polarization (see \cite[\S2.2.3]{PV})\footnote{The volume is independent of the choice of weak polarization, see \cite[Lemma 2.2.2]{PV}.}.
\end{conj2}
\begin{prop2}Assuming the Beilinson conjecture, Conjecture \ref{Period} is equivalent to the isometry statement in Conjecture \ref{Metric} for top and bottom degrees,
\[\wedge^{\mathrm{top}}\cM^{*}\otimes\cH^{i_{\min}}\cong \cH^{i_{\max}}.\] In particular, if $i_{\max}=i_{\min}+1$, Conjecture \ref{Metric} and Conjecture \ref{Period} are equivalent.
\end{prop2}
\begin{proof}
The key is to compare the top and the bottom degrees of the graded vector spaces and relate them with the volumes, which is independent of the choice of metric. Namely, we need to prove the analogue of \cite[Lemma 2.2.2]{PV},
\[\vol_{S} H^{1}_{M}((\Ad\Pi)_{\cO_{E}},\ov{\bQ}(1))\sim_{\ov{\bQ}^{\times}}\frac{L^{*}(0,\Ad\Pi)}{\vol_{S} F^{1}H_{\dR}(\Ad\Pi)},\]for $S$ the Hermitian inner product induced by a weak polarization, and the statement will follow after applying functional equation. This now follows from the Beilinson's conjecture over a general number fields, as in \cite[\S6]{Ramakrishnan}. Namely, there is a fundamental exact sequence, \cite[(6.4.2)]{Ramakrishnan},
\[0\rar F^{1}H_{\dR}((\Ad\Pi)_{\bC})\rar H_{B}^{0}((\Ad\Pi)_{\bR},\bC)\rar H_{\sD}^{1}((\Ad\Pi)_{\bR},\bC(1))\rar0\tag{$\star$}\label{FES},\]where\[H_{B}^{0}((\Ad\Pi)_{\bR},\bC)=\prod_{w\text{ complex places of $E$}}H_{B}^{0}(M\times_{E,w}\bC,\bC)\times\prod_{w\text{ real places of $E$}}H_{B}^{0}(M\times_{E,w}\bR,\bC),\]is similarly defined as the ``real Deligne cohomology''. The Beilinson's conjecture over a general number field says that the determinant of the fundamental exact sequence (\ref{FES}) has incompatible $\ov{\bQ}$-rational structures, and are off precisely by $L^{*}(0,\Ad\Pi)$ (see \cite[\S6.4]{Ramakrishnan}):
\[\det(H_{B}^{0}((\Ad\Pi)_{\bR},\ov{\bQ}))L^{*}(0,\Ad\Pi)\sim_{\ov{\bQ}^{\times}}\det F^{1}H_{\dR}((\Ad\Pi)_{\ov{\bQ}})\cdot\det(H_{M}^{1}((\Ad\Pi)_{\cO_{E}},\ov{\bQ}(1))).\]Regarding this as an equality inside the determinant of the fundamental exact sequence (\ref{FES}), computing the volumes would give the desired statement.
\end{proof}We will later prove Conjecture \ref{Metric} in special cases where Conjecture \ref{Period} is equivalent to Conjecture \ref{Metric} (namely, the appearances of Hecke eigensystem span over two degrees).
\begin{exam2}[(Sanity check: $\SL_{2}$)]\label{SL2}
The simplest example is the case of $G=\SL_{2,\bQ}$, where the conjecture is about \emph{weight one elliptic modular forms}. Let $f\in S_{1}(\Gamma)$ be a weight one cuspidal new eigenform with Fourier coefficients in $\ov{\bQ}$, generating an automorphic representation $\Pi$. Then, the complex conjugate $\ov{f}\in\ov{\Pi}$, which also satisfies $\langle\ov{f},\ov{f}\rangle_{P}=\langle f,f\rangle_{P}$. Thus, we are led to the following question: for what $c\in\bC^{\times}$ does $c\ov{f}$ define a $\ov{\bQ}$-coherent cohomology class in $H^{1}(X(\Gamma),\omega)$? After calculating the archimedean $L$-factors (which is elementary), Conjecture \ref{Period} says that
\[|c|^{-2}\sim_{\ov{\bQ}^{\times}\cap\bR}\left(\pi^{-2}\frac{L(1,\Pi,\Ad)}{\vol F^{1}H_{\dR}(\Ad\Pi)}\right)^{2}.\]Note that $\Pi$ is a pure weight $0$ motive, and so is $\Ad\Pi$; therefore, $F^{1}H_{\dR}(\Ad\Pi)=0$, which means that
\[|c|\sim_{\ov{\bQ}^{\times}\cap\bR}\frac{\pi^{2}}{L(1,\Pi,\Ad)}.\] It is well-known that $L(1,\Pi,\Ad)\pi^{-2}\sim_{\ov{\bQ}^{\times}}\langle f,f\rangle_{P}$. Thus, the conjecture says 
\[|c|\sim_{\ov{\bQ}^{\times}\cap\bR}\frac{\pi}{\langle f,f\rangle_{P}}.\]Let $c\ov{f}d\ov{z}$ define a Dolbeault cohomology class $f^{\vee}\in H^{1}(X(\Gamma),\omega)$, which is, by definition, defined over $\ov{\bQ}$. By Serre duality, $\langle f,f^{\vee}\rangle_{S}\in\ov{\bQ}$, where $\langle-,-\rangle_{S}$ denotes the Serre duality. 

On the other hand, the Serre duality in this case coincides with Petersson inner product scaled by the factor of $\frac{1}{2\pi i}$ (e.g. \cite[p. 22]{Deligne}), namely $2\pi i\langle f,f^{\vee}\rangle_{S}=\langle f,\ov{c}f\rangle_{P}=\ov{c}\langle f,f\rangle_{P}$. Thus, means that $c\sim_{\ov{\bQ}^{\times}}\frac{1}{\langle f,f\rangle_{P}}$,\footnote{Indeed, $\langle f,f\rangle_{P}$ is a real number.} which is consistent with our conjecture. That these facts can be realized as an instance of motivic action conjecture was realized in \cite{HV} and was explicitly spelled out in \cite{Horawa}.
\end{exam2}
\subsection{An approach towards the motivic action conjecture}
Unlike \cite{PV} or the case of $\SL_{2}$, the top and bottom degrees are most of the time not complementary (namely, $i_{\Pi_{\min}}+i_{\Pi_{\max}}\ne\dim_{\bC}X$), so the conjecture has nothing to do with any form of duality. We will nevertheless prove a form of the conjecture by relating this with \emph{cohomological period integrals}. These are integral representations of certain $L$-functions that \begin{itemize}\item apply to automorphic forms appearing in higher coherent cohomology,
\item and admit an interpretation as cup product pairing in coherent cohomology.
\end{itemize}
This is useful in verifying our conjectures as the coherent cohomological cup product can \emph{detect rationality of higher coherent cohomology classes}.

We will show that a strategy similar to \cite[\S7]{PV} can also show that our period conjecture, Conjecture \ref{Period}, is true in certain cases.
\begin{thm2}\label{MainEvidence}
Let $G$ be either $\Sp_{4}$ or $\SU(2,1)$, and let $\Pi$ be a globally generic cuspidal automorphic representation of $G(\bA_{\bQ})$ satisfying Assumption \ref{Assumption}. Assume the working hypothesis on periods, Assumption \ref{PeriodAssumption}. Then, Conjecture \ref{Period} is true, up to the factor of an archimedean zeta integral (see Remarks \ref{Sp4Whittaker}, \ref{SU21Whittaker}).
\end{thm2} In the two cases, we will compare periods of \emph{holomorphic LDS} appearing in $H^{0}$ and \emph{generic LDS} appearing $H^{1}$. The tools that we will use are summarized in the following table.
\renewcommand{\arraystretch}{1.5}
\begin{center}
\resizebox{\textwidth}{!}{\begin{tabular}{|c|c|c|}
\hline
 LDS type & Detecting rationality & Relation with Petersson norm\\
 \hline
 Holomorphic ($H^{0}$)& Rational Fourier coefficient & {\renewcommand{\arraystretch}{1.2}\begin{tabular}{@{}c@{}}Doubling method, Refined GGP conjectures \end{tabular}}\\
 \hline
 Generic ($H^{>0}$) & Cohomological period integrals & Lapid--Mao conjecture \\
 \hline
 \end{tabular}}
 \end{center}
 A slightly more detailed outline is as follows. We would need to know how $\ov{\bQ}$-algebraicity of coherent cohomology classes in $H^{0}$ and $H^{1}$ is related to Petersson norms. For $H^{0}$, the classes are represented by holomorphic automorphic forms, where their algebraicity is detectable by Fourier coefficients. In the language of periods, these are related to Bessel or Fourier--Jacobi periods, which are related to Petersson norms via the refined Gan--Gross--Prasad conjectures. For $H^{1}$, the classes are represented by generic automorphic forms\footnote{In both cases of our concern, the infinity type corresponding to $H^{1}$ belongs to a generic (L)DS, which is a numerical coincidence that only happens in certain special examples.}. The Petersson norms of corresponding automorphic forms are related to the \emph{Whittaker periods} via the \emph{Lapid--Mao conjecture} \cite{LapidMao}. Due to its simple statement, we recall the conjecture here:\begin{conj2}[({\cite{LapidMao}})]\label{LMConj}Let $\Pi$ be a globally generic representation, satisfying (\ref{NEW}) of Assumption \ref{Assumption}. Let $f=\otimes_{p}f_{p}\in\Pi$ be a newform. Then,
\[\frac{\langle f,f\rangle}{|W(1)|^{2}}\sim_{\bQ^{\times}}\frac{L(1,\Pi,\Ad)}{\Delta_{G}(1)|W_{\infty}(1)|^{2}},\]where $\langle,\rangle$ is the $L^{2}$-norm (as opposed to the Petersson norm), $\Delta_{G}(s)$ is the $L$-function of the dual to the Artin motive attached to $G$ as defined in \cite{Gross}, and $W$ ($W_{\infty}$, respectively) is the Whittaker function of $f$ ($f_{\infty}$, respectively).
\end{conj2}The Whittaker periods are then related to $\ov{\bQ}$-algebraicity of coherent cohomology classes via cohomological period integrals. On one hand, the period integral has automorphic interpretation, which connects to Whittaker periods. On the other hand, the period integral has cohomological interpretation, so that in particular it descends to $\ov{\bQ}$, hence detects $\ov{\bQ}$-algebraicity. 

Finally, as the conjectures are formulated using invariants coming from the motivic formalism, we would be working with the corresponding motives and compute motivic invariants (dubbed ``Hodge-linear algebra''). We will work with the adjoint motive, as in Conjecture \ref{AdjointMotive}, as well as the motive of the given automorphic representation\footnote{The construction of those motives is a subtle matter, as we work with limits of discrete series. Indeed, the corresponding Galois representations have been constructed, but only by using congruences.}:
\begin{conj2}\label{ActualMotive}For $\Pi$ as in Assumption \ref{Assumption}, there exists a motive $M_{\Pi}$ that is uniquely characterized by \cite[\S4.3.3]{Clozel}.
\end{conj2}
\begin{rmk2}The adjoint motive of Conjecture \ref{AdjointMotive} is the motive associated to $L(s,\Pi,\Ad)$ in the above sense.
\end{rmk2}
Thus, the ``working hypothesis'' is as follows.
\begin{ass2}[(Working hypothesis)]\label{PeriodAssumption}
We assume the following conjectures\footnote{Even though the refined Gan--Gross--Prasad conjectures are being mentioned throughout the paper, they are not required as an assumption, because we have an unconditional alternative result using the doubling method (and they are equivalent if one assumes Beilinson's conjectures anyways).}. There are numerous instances of these conjectures being verified, and we do not attempt to list them here.
\begin{enumerate}
\item \emph{Beilinson's conjecture for Chow motives}, Conjecture 2.1.1 of \cite{PV}.
\item \emph{Lapid--Mao conjecture}, Conjecture \ref{LMConj}.
\item \emph{Existence of motives}, Conjecture \ref{AdjointMotive} and Conjecture \ref{ActualMotive}.
\end{enumerate}
\end{ass2}
\section{Evidence I: The case of $\Sp_{4}$}
\label{Sp4EvidenceSection}
We first provide the evidence for the motivic action conjecture (Conjecture \ref{Metric}) for the case of certain irregular automorphic forms on $G(\bA)=\Sp_{4}(\bA)$. In this case, we will use an integral reprsentation of the spinor $L$-function by Novodvorsky (\cite{Nov}), whose coherent cohomological interpretation was given by \cite{LPSZ}. 

Let $\Pi$ be a globally generic cuspidal automorphic representation of $G=\GSp_{4}(\bA_{\bQ})$, namely that the Whittaker transform (see \S\ref{Notation}, Notation) defines a realization of $\Pi$ as a space of functions on $G(\bQ)\bs G(\bA)$ satisfying a transformation property under the $N(\bA)$-action. We choose unramified vectors $\varphi_{v}^{0}\in \Pi_{v}$ for $v$ finite with $\Pi_{v}$ unramified such that, if $\psi_{v}$ is unramified, $W_{\varphi_{v}^{0}}(1)=1$. Suppose that $\Pi$ satisfies Assumption \ref{Assumption}.
\begin{defn2}Let $M\supset F_{\Pi}$ be a number field. A $\psi$-Whittaker function $W$ on $G$ is called to be \emph{defined over $M$} if it takes values in $M(\mu_{\infty})$ and satisfies 
\[\sigma(W(g))=W(w(\kappa(\sigma))g),\]for all $g\in G(\bA_{f})$ and $\sigma\in\Gal(\ov{\bQ}/M)$, where $w(x)=\diag(x^{3},x^{2},x,1)$ and $\kappa:\Gal(\ov{\bQ}/\bQ)\rar\wh{\bZ}^{\times}$ is the cyclotomic character.
\end{defn2}
\begin{rmk2}\label{RKW}
\hfill\begin{enumerate}
\item The convention is made so that, if $\psi$ is unramified, the unramified $\psi$-Whittaker function $W$ with $W(1)=1$ is defined over $F_{\Pi}$. As $L(s,\Pi,\Ad)$ is regular at $s=1$, $W(1)$ is nonzero by \cite[\S3.1]{LapidMao}. Thus, $W$ is defined over $M$ if and only if $W(1)\in M^{\times}$. 
\item Indeed, our conditions on $G$ and $\Pi$ imply that $F_{\Pi}$ is a number field (see Notation). Thus, the space of $\psi$-Whittaker functions, which naturally has an action of $\Aut(\bC/\bQ)$ as in \cite[\S4.1]{GHL}, is fixed by $\Aut(\bC/F_{\Pi})$. Thus, there exists a nontrivial $\psi$-Whittaker function defined over $M$.
\end{enumerate}
\end{rmk2}
\subsection{Whittaker periods via cohomological period integrals}The purpose of this subsection is to prove the following
\begin{thm2}\label{Sp4W}Let $\varphi\in\Pi^{\new}=\Pi_{f}^{\new}\otimes\Pi_{\infty}^{\new}$. Let $F$ be an imaginary quadratic field such that the special Bessel period of $\Pi$ for $F$ is not identically zero (see (\ref{Sp4A}) of Assumption \ref{Assumption}). If $[\varphi]$ (see Definition \ref{CohclassDef}) is defined over a number field $F'\supset FF_{\Pi}$, then 
${\Lambda\left(1/2,\Pi\right)\Lambda\left(1/2,\Pi\otimes\chi_{F}\right)}W_{\varphi}$ is a nontrivial $\psi$-Whittaker function defined over $F'$, where $\chi_{F}$ is the quadratic character associated to $F$.
\end{thm2}
\begin{proof}
Let $C\in\bC^{\times}$ be such that $\frac{W_{\varphi}}{C}$ is defined over $F'$, which exists by Remark \ref{RKW}(2). We would like to show that \[C\Lambda(1/2,\Pi)\Lambda(1/2,\Pi\otimes\chi_{F})\in F'{}^{\times}.\]We first prove that the above quantity is nonzero. This is because, by \cite[(1.26)]{BochererFurusawa}, the nonvanishing of $L(1/2,\Pi)L(1/2,\Pi\otimes\chi_{F})$ is equivalent to nonvanishing of the special Bessel period $B(f_{\Pi},F)$ in (\ref{Sp4A}) of Assumption \ref{Assumption}; that the special Bessel period defined in this paper is the Bessel period (up to explicit nonzero scalar) in \cite{BochererFurusawa} follows from the calculation of \cite[Proposition 3.5]{BochererRefined}.

We now consider Novodvorsky's integral representation of spinor $L$-function. This is, roughly speaking, the period integral of $\varphi$ times an Eisenstein series over an embedded product of two modular curves. More precisely, there is an embedding of $H=\GL_{2}\times_{\GL_{1}}\GL_{2}$ into $G$, which is most naturally thought as $\SO(2,2)\hrar\SO(3,2)$. Let $B\subset H$ be the upper triangular Borel. For $\Phi_{1},\Phi_{2}:\bA^{2}\rar\bC$ and $\chi_{1},\chi_{2}$ a unitary Gr\"ossencharcter, we define an Eisenstein series on $H$ with respect to $B$, $E(h,\chi_{i},\Phi_{i},s_{i})$. The following are well-known (e.g. \cite[Proposition 7.3]{LPSZ}).
\begin{prop2}\label{Shim}Suppose we are taking the ``weight $k$-section'' for $\Phi_{i,\infty}$.\begin{enumerate}
\item For $-\frac{k}{2}+1\le s_{1},s_{2}\le\frac{k}{2}$ half-integral, equivalent modulo $1$ to $\frac{k}{2}$, this defines a nearly holomorphic form on $H$, which can be thought as an $H^{0}$ class of some automorphic vector bundle over a Shimura variety for $H$.
\item If $\chi_{i}$'s and $\Phi_{i}$'s are valued in a number field $F$, then the $H^{0}$ class is defined over $F$.
\end{enumerate}
\end{prop2}
Let $(\lambda_{1},\lambda_{2})$ be the Harish--Chandra character of the (necessarily generic) (L)DS $\Pi_{\infty}$. Let $\Gamma_{H}=H\cap\Gamma$, and let $i:X_{H}(\Gamma_{H})\hrar X_{G}(\Gamma)$ be the closed embedding of certain toroidal compactifications of the corresponding Shimura varieties\footnote{That the closed immertion of open Shimura varieties extends to a closed immersion of toroidal compactifications with respect to certain refinements of polyhedral cone decomposition is achieved by \cite{LanClosed}, but the situation is simplier in this case, because $X_{H}$, being a productof modular curves, is unique.}. Then, Novodvorsky's integral representation can be understood via a cup product pairing \[H^{2}(X_{G},V)\otimes_{\bC}H^{0}(X_{H},W)\xrar{(\id,i_{*})}H^{2}(X_{G},V)\otimes_{\bC}H^{1}(X_{G},W')\xrar{\cup} H^{3}(X_{G},V\otimes W')\xrar{S} H^{0}(X_{G},\cO)=\bC,\]where $\cup$ is the cohomological cup product and $S$ is the Serre duality pairing, induced from a morphism of algebraic $\ov{\bQ}$-representations of $K_{\infty}$, $V\otimes W'\rar\fg^{-1,1}$ (namely, the pairing is normalized such that, on the level of representations of $K_{\infty}$, the $\ov{\bQ}$-structures are compatible). Here, $W$ and $W'$ are certain automorphic vector bundles over $X_{H}$ and $X_{G}$, respectively, corresponding to $[V_{H}(\lambda_{1}-\lambda_{2}-1,0)]\otimes\omega_{H}(1,1)$ and $[\wt{L}_{1}]$, if we use the notation of \cite[\S6]{LPSZ}. The reinterpretation, done in \cite[\S7.4]{LPSZ}, of Novodvosky's integral asserts that, given an imaginary quadratic field $F$,  
there is a coherent cohomology class $[E]\in H^{0}(X_{H}(\Gamma_{H}),W)$ defined over $FF_{\Pi_{f}}$, which corresponds to a nearly-holomorphic Eisenstein series $E$ under the Hodge splitting of \cite[\S6.3]{LPSZ}, 
such that \[\langle[\varphi],[E]\rangle=C\Lambda\left(1/2,\Pi\right)\Lambda\left(1/2,\Pi\otimes\chi_{F}\right),\] where $\chi_{F}$ is the Hecke character corresponding to $F$.
As $\langle,\rangle$ is defined over $\bQ$, the left hand side is in $F'$, which is what we are looking for.
\end{proof}

\subsection{Hodge-linear algebra}\label{Sp4Hodge}
We now put the relevant Hodge-linear algebra that proves the case of $\Sp_{4}$. For the sake of simplicity, we assume we work with a parallel weight $(2,2)$ Siegel modular form, or Harish--Chandra parameter $(1,0)$, although the Hodge-linear algebra calculation stays the same for general weights. Our goal is to convert, using elementary linear algebra, $\vol F^{1}H_{\dR}(\Ad M)$ into an expression that involves Deligne's periods $c^{+}(M),c^{-}(M),\delta(M)$ whose definitions will be recalled later. We will then be able to express $\vol F^{1}H_{\dR}(\Ad M)$ with $L$-values, using Deligne's conjectures. We will prove
\begin{prop2}\label{Sp4HodgeProp}For a motive $M$ associated to $\Pi$ (in the sense of Conjecture \ref{ActualMotive}), we have
\[\vol F^{1}H_{\dR}(\Ad M)\sim_{\ov{\bQ}^{\times}}\frac{\sqrt{c^{+}(M)c^{-}(M)}^{3}}{\delta(M)^{3/2}}.\]
\end{prop2}
\begin{proof}
The motive $M$ of $\Pi$ should be of rank $4$ and weight $1$, with the Hodge decomposition
\[H_{B}(M)\otimes_{\bQ}\bC=H^{1,0}(M)\oplus H^{0,1}(M),\quad \dim_{\bC}H^{1,0}(M)=\dim_{\bC}H^{0,1}(M)=2,\]which is of the type of the Hodge structure defined by the corresponding archimedean Langlands parameter.
In this case, \[\delta(M)\text{ }\left(\text{$c^{\pm}(M)$, resp.}\right)\text{ }\in\bC^{\times}/\bQ^{\times},\]is the determinant of the comparison map \[\toolong{H_{B}(M)\otimes\bC\riso H_{\dR}(M)\otimes\bC\text{ }\left(\text{$H_{B}(M)^{+}\otimes\bC\rar H_{B}(M)\otimes\bC\riso H_{\dR}(M)\otimes\bC\thrar (H_{\dR}(M)/F^{1}H_{\dR}(M))\otimes\bC$, resp.}\right),}\]with respect to the bases coming from the underlying $\bQ$-structures on both sides.

Let $e_{1}^{+},e_{2}^{+}$ be a $\bQ$-basis of $H_{B}(M)^{+}$ and $e_{1}^{-},e_{2}^{-}$ be a $\bQ$-basis of $H_{B}(M)^{-}$. Let $f_{1},f_{2}\in F^{1}H_{\dR}(M)$ be a $\bQ$-basis, and $g_{1},g_{2}\in H_{\dR}(M)/F^{1}H_{\dR}(M)$ be a $\bQ$-basis, and $\wt{g_{1}},\wt{g_{2}}\in H_{\dR}(M)$ be lifts of $g_{1},g_{2}$. Given two $\bC$-bases of $H_{B}(M)\otimes\bC\cong H_{\dR}(M)\otimes\bC$, we can write an expression
\[\begin{pmat}e_{1}^{+}& e_{2}^{+}& e_{1}^{-}& e_{2}^{-}\end{pmat}=\begin{pmat}f_{1}& f_{2} & \wt{g_{1}}& \wt{g_{2}}\end{pmat}\begin{pmat}A&B\\C&D\end{pmat},\]for $A,B,C,D\in M_{2}(\bC)$. Note that by definition $\delta(M)=\det\begin{psmat}A&B\\C&D\end{psmat}$.

We have canonical isomorphisms
\[H^{1,0}(M)\cong F^{1}H_{\dR}(M)\otimes\bC,\quad H^{0,1}(M)\cong\frac{H_{\dR}(M)}{F^{1}H_{\dR}(M)}\otimes\bC.\]Let $f_{1,B},f_{2,B}\in H^{1,0}(M)$ and $g_{1,B},g_{2,B}\in H^{0,1}(M)$ be the images of $f_{1},f_{2},g_{1},g_{2}$ under the above canonical isomorphisms. Then $f_{i,B}$ and $f_{i}$ coincide as elements of $H_{B}(M)\otimes\bC\cong H_{\dR}(M)\otimes\bC$, whereas $\wt{g_{i}}-g_{i,B}\in H^{1,0}(M)$. So, \[\begin{pmat}f_{1}&f_{2}&\wt{g_{1}}&\wt{g_{2}}\end{pmat}=\begin{pmat}f_{1,B}&f_{2,B}&g_{1,B}&g_{2,B}\end{pmat}\begin{pmat}1_{2}&M\\0_{2}&1_{2}\end{pmat}.\]In particular, if we write\[\begin{pmat}e_{1}^{+}& e_{2}^{+}& e_{1}^{-}& e_{2}^{-}\end{pmat}=\begin{pmat}f_{1,B}& f_{2,B} & g_{1,B}& g_{2,B}\end{pmat}\begin{pmat}A'&B'\\C'&D'\end{pmat},\]then $\det\begin{psmat}A&B\\C&D\end{psmat}=\det\begin{psmat}A'&B'\\C'&D'\end{psmat}$. 

 Also, there must be relations
\[c_{B}(f_{1,B})=ag_{1,B}+bg_{2,B},\]
\[c_{B}(f_{2,B})=cg_{1,B}+dg_{2,B}.\]Using that $F_{\infty}e_{i}^{+}=e_{i}^{+}$, $F_{\infty}e_{i}^{-}=-e_{i}^{-}$ and that $F_{\infty}(f_{i,B})=c_{B}(f_{i,B})$ and $F_{\infty}(g_{i,B})=c_{B}(g_{i,B})$ (as in \cite[Lemma, \S8.2.1]{PV}), one has
\begin{eqnarray*}c_{B}(A'_{1i}f_{1,B}+A'_{2i}f_{2,B})&=&C'_{1i}g_{1,B}+C'_{2i}g_{2,B},\\c_{B}(B'_{1i}f_{1,B}+B'_{2i}f_{2,B})&=&-(D'_{1i}g_{1,B}+D'_{2i}g_{2,B}),\end{eqnarray*}for $i=1,2$. This can be packaged into $C'=\begin{psmat}a&c\\b&d\end{psmat}A'$, $D'=-\begin{psmat}a&c\\b&d\end{psmat}B'$. So
\[\delta(M)=\det\begin{pmat}A'&B'\\\begin{psmat}a&c\\b&d\end{psmat}A'&-\begin{psmat}a&c\\b&d\end{psmat}B'\end{pmat}=\det\begin{pmat}A'&B'\\0_{2}&-2\begin{psmat}a&c\\b&d\end{psmat}B'\end{pmat}=4\det\begin{pmat}a&c\\b&d\end{pmat}\det A'\det B'.\]Note on the other hand that $c^{+}(M)=\det C$, $c^{-}(M)=\det D$. As $C'=C$, $D'=D$, we see that
\[c^{+}(M)c^{-}(M)\sim_{\ov{\bQ}^{\times}}\delta(M)\det\begin{pmat}a&c\\b&d\end{pmat}.\]

Note that, as $M$ is self-dual, $M^{\vee}\cong M(1)$, which implies $\Ad M=(\Ad M)^{*}=(\Sym^{2}M)(1)$. We are also led to calculate $\vol F^{1}H_{\dR}(\Ad M)$ (we know it does not depend on the choice of weak polarization up to $\ov{\bQ}^{\times}$-ambiguity). We know that $(\vol F^{1}H_{\dR}(\Ad M))^{2}\sim_{\ov{\bQ}^{\times}}\lambda$, where $\varphi(c_{B}(v^{+}))=\lambda v^{-}$. Here, $v^{+},v^{-}$ are $\bQ$-bases vectors for $\det F^{1}H_{\dR}(\Ad M)$ and $\det\left(\frac{H_{\dR}(\Ad M)}{F^{0}H_{\dR}(\Ad M)}\right)$, respectively, and $\varphi:\wedge^{\dim F^{1}H_{\dR}(\Ad M)}(H_{\dR}(M)\otimes\bC)\rar \det\left(\frac{H_{\dR}(\Ad M)}{F^{0}H_{\dR}(\Ad M)}\right)$ is the natural projection. As $\Ad M=(\Sym^{2} M)(1)$, we can take $f_{1}^{2},f_{1}f_{2},f_{2}^{2}$ and $g_{1}^{2},g_{1}g_{2},g_{2}^{2}$ as $\bQ$-bases of $F^{1}H_{\dR}(\Ad M)$ and $\frac{H_{\dR}(\Ad M)}{F^{0}H_{\dR}(\Ad M)}$, respectively. Now from the known relations,
\begin{eqnarray*}c_{B}(f_{1}^{2})&\equiv&a^{2}g_{1}^{2}+2abg_{1}g_{2}+b^{2}g_{2}^{2},\\c_{B}(f_{1}f_{2})&\equiv&acg_{1}^{2}+(ad+bc)g_{1}g_{2}+bdg_{2}^{2},\\c_{B}(g_{1}g_{2})&\equiv&c^{2}g_{1}^{2}+2cdg_{1}g_{2}+d^{2}g_{2}^{2},\end{eqnarray*}where $\equiv$ is mod $F^{0}H_{\dR}(\Ad M)\otimes\bC$. So \[\lambda=\det\begin{pmat}a^{2}&2ab&b^{2}\\ac&ad+bc&bd\\c^{2}&2cd&d^{2}\end{pmat}=(ad-bc)^{3}.\]Therefore, we obtain
\[\vol F^{1}H_{\dR}(\Ad M)\sim_{\ov{\bQ}^{\times}}\frac{\sqrt{c^{+}(M)c^{-}(M)}^{3}}{\delta(M)^{3/2}}.\]
\end{proof}\subsection{Completion of the proof}
\begin{proof}[Proof of Theorem \ref{MainEvidence} for $\Sp_{4}$]
Now we apply the relevant period conjectures for this case. Let $f_{\hol}$, $f_{\gen}$ be newforms (see Assumption \ref{Assumption}) in $\Pi_{f}\otimes\Pi_{\hol}$, $\Pi_{f}\otimes\Pi_{\gen}$, respectively, where $\Pi_{\hol}$ is the corresponding holomorphic NLDS in the $L$-packet of $\Pi_{\gen}$. Let us assume that $[f_{\hol}]$ and $[f_{\gen}]$ are defined over $\ov{\bQ}$. Let $F$ be the imaginary quadratic field as in Theorem \ref{Sp4W}. Then, a theorem of Furusawa--Morimoto\footnote{See \cite[Theorem 1]{BochererFurusawa} for the case of discrete series; the case of limit of discrete series is also an upcoming work of them.} \footnote{This is a specific case of the refined Gan--Gross--Prasad conjecture for Bessel periods, \cite[Conjecture 2.5]{Liu}. This special case is also sometimes called \emph{B\"ocherer's conjecture}.} implies that\[\frac{\langle f_{\hol},f_{\hol}\rangle}{|B(f_{\hol},F)|^{2}}\sim_{\ov{\bQ}^{\times}}\pi^{-6}\frac{L(1,\Pi,\Ad)}{L(1/2,\Pi)L(1/2,\Pi\otimes{\chi_{F}})}.\tag{A}\]
Since $[f_{\hol}]$(see Definition \ref{CohclassDef}) is defined over $\ov{\bQ}$, we know that $B(f_{\hol},F)$, a $\bQ$-linear combination of Fourier coefficients of $f_{\hol}$, is in $\ov{\bQ}$. Thus, in the setting of Theorem \ref{MainEvidence}, $\langle f_{\hol},f_{\hol}\rangle\sim_{\ov{\bQ}^{\times}}\pi^{-6}\frac{L(1,\Pi,\Ad)}{L(1/2,\Pi)L(1/2,\Pi\otimes\chi_{F})}$. 

We now relate the $\ov{\bQ}$-rationality of $f_{\gen}$ with $\langle f_{\gen},f_{\gen}\rangle$. Note that Theorem \ref{Sp4W} is about the relationship between rationality of Whittaker functions and that of coherent cohomology classes, for the anti-generic LDS (namely, those appearing in $H^{2}$ of coherent cohomology). Thus, as $[f_{\gen}]$ is defined over $\ov{\bQ}$, by Serre duality,
\[\label{AA}\langle f_{\gen},f_{\gen}\rangle_{P}=(2\pi i)^{3}\langle[f_{\gen}],[\ov{f}_{\gen}]\rangle_{\coh},\tag{B}\]where $\langle-,-\rangle_{\coh}$ is the cohomological cup product defined over $\ov{\bQ}$ normalized so that it is induced from the $\ov{\bQ}$-morphism of algebraic $K_{\infty}$-representations $W\otimes\Hom(W,\fg^{-1,1})\rar\fg^{-1,1}$. If we let $C\in\bC^{\times}$ be such that $C[\ov{f}_{\gen}]$ is defined over $\ov{\bQ}$, we see that the RHS is $\sim_{\ov{\bQ}^{\times}}\pi^{3}C^{-1}$. On the other hand, by Theorem \ref{Sp4W}, \[\Lambda(1/2,\Pi)\Lambda(1/2,\Pi\otimes\chi_{F})W_{C\ov{f}_{\gen}}=C\Lambda(1/2,\Pi)\Lambda(1/2,\Pi\otimes\chi_{F})W_{\ov{f}_{\gen}},\]is defined over $\ov{\bQ}$. We now invoke the Lapid--Mao conjecture, Conjecture \ref{LMConj}:
 \[\frac{\langle f_{\gen},f_{\gen}\rangle}{|W(1)|^{2}}\sim_{\ov{\bQ}^{\times}}\pi^{9}\cdot\frac{L(1,\Pi,\Ad)}{|W_{\infty}(1)|^{2}}.\]By Remark \ref{RKW}(1), $W_{\ov{f}_{\gen}}(1)\ne0$, so that \[\label{BB}C\Lambda(1/2,\Pi)\Lambda(1/2,\Pi\otimes\chi_{F})W_{\ov{f}_{\gen}}(1)\in\ov{\bQ}^{\times}.\tag{C}\]Since $W_{\ov{f}_{\gen}}(1)=\ov{W_{f_{\gen}}(1)}$, we have
\[\pi^{9}\frac{L(1,\Pi,\Ad)}{|W_{\infty}(1)|^{2}}|W_{f_{\gen}}(1)|^{2}\overtext{Lapid--Mao}{\sim_{\ov{\bQ}^{\times}}}\langle f_{\gen},f_{\gen}\rangle\overtext{(B)\quad }{\sim_{\ov{\bQ}^{\times}}}\pi^{3}C^{-1}\overtext{(C)\quad }{\sim_{\ov{\bQ}^{\times}}}\pi^{3}\ov{W_{f_{\gen}}(1)}\Lambda(\frac{1}{2},\Pi)\Lambda(\frac{1}{2},\Pi\otimes\chi_{F}),\]
or
\[W_{f_{\gen}}(1)\sim_{\ov{\bQ}^{\times}}\pi^{-6}\frac{|W_{\infty}(1)|^{2}\Lambda(\frac{1}{2},\Pi)\Lambda(\frac{1}{2},\Pi\otimes\chi_{F})}{L(1,\Pi,\Ad)}\sim_{\ov{\bQ}^{\times}}\pi^{-10}\frac{|W_{\infty}(1)|^{2}L(\frac{1}{2},\Pi)L(\frac{1}{2},\Pi\otimes\chi_{F})}{L(1,\Pi,\Ad)}.\tag{D}\]
We need to compute $\frac{\langle f_{\hol},f_{\hol}\rangle_{P}}{\langle f_{\gen},f_{\gen}\rangle_{P}}$, which can be now seen as follows.
\begin{align*}\frac{\langle f_{\hol},f_{\hol}\rangle_{P}}{\langle f_{\gen},f_{\gen}\rangle_{P}}&\overtext{Lapid--Mao + (A)}{\sim_{\ov{\bQ}^{\times}}}\frac{\pi^{-6}\frac{L(1,\Pi,\Ad)}{L(1/2,\Pi)L(1/2,\Pi\otimes\chi_{F})}}{\pi^{9}|W_{f_{\gen}}(1)|^{2}\frac{L(1,\Pi,\Ad)}{|W_{\infty}(1)|^{2}}}
\\&=\pi^{-15}\frac{|W_{\infty}(1)|^{2}}{|W_{f_{\gen}}(1)|^{2}L(1/2,\Pi)L(1/2,\Pi\otimes\chi_{F})}\\&\overtext{(D)}{\sim_{\ov{\bQ}^{\times}}}\pi^{-15}\frac{|W_{\infty}(1)|^{2}}{\pi^{-20}\frac{|W_{\infty}(1)|^{4}L(1/2,\Pi)^{3}L(1/2,\Pi\otimes\chi_{F})^{3}}{L(1,\Pi,\Ad)^{2}}}\\&\sim_{\ov{\bQ}^{\times}}\frac{\pi^{9}}{|W_{\infty}(1)|^{2}}\left(\frac{L_{\infty}(1,\Pi,\Ad)}{L_{\infty}(0,\Pi,\Ad)}\cdot\frac{L(1,\Pi,\Ad)}{\frac{\sqrt{L(1/2,\Pi)L(1/2,\Pi\otimes\chi_{\bQ(i)})}^{3}}{\pi^{3}}}\right)^{2},\tag{E}\label{FINALSP4}
\end{align*}
which involves elementary calculation of archimedean $L$-factors. 

On the other hand, by Deligne's conjectures, Proposition \ref{Sp4HodgeProp} implies that
\[\vol F^{1}H_{\dR}(\Ad M)\sim_{\ov{\bQ}^{\times}}\frac{\sqrt{L(1/2,\Pi)L(1/2,\Pi\otimes\chi_{\bQ(i)})}^{3}}{\pi^{3}}.\]
Thus, the parenthized term in (\ref{FINALSP4}) is precisely the RHS of Conjecture \ref{Period}, which finishes the proof.
\end{proof}
\begin{rmk2}\label{Sp4Whittaker}
Unfortunately, for now, the author has been unable to calculate $W_{\infty}(1)$. It is however very believable that the archimedean zeta integral against a preferred, nice test vector is equal to an archimedean local $L$-factor, which is $\sim_{\ov{\bQ}^{\times}}$ half-integral powers of $\pi$. In our case, there is even an explicit integral expression of $W_{\infty}(1)$, as written in \cite{GenericBocherer}:
\[W_{\infty}(1)=16e^{-2\pi}\pi^{\frac{7}{2}}\int_{c-i\infty}^{c+i\infty}\pi^{-2s}\Gamma(s+\frac{1}{2})^{2}U(s+\frac{1}{2},1,4\pi)\Gamma(s)\frac{ds}{2\pi i},\]where $U(a,b,z)=\frac{1}{\Gamma(a)}\int_{0}^{\infty}e^{-zt}t^{a-1}(1+t)^{b-a-1}dt$ is the confluent hypergeometric function of the second kind.
\end{rmk2}

\section{Evidence II: The case of $\SU(2,1)$}\label{SU21EvidenceSection}In this section, we provide the evidence for the Period conjecture (Conjecture \ref{Period}) for the case of certain irregular automorphic forms on $G(\bA)=\SU(2,1)(\bA)$. In this case, we will use an integral reprsentation of the base-change $L$-function by Gelbart and Piatetski-Shapiro (\cite{GPS}, completed by \cite{KO}), whose coherent cohomological interpretation was given by \cite{Oh}.

Let $\Pi$ be a globally generic cuspidal automorphic representation of $G=\SU(2,1)(\bA_{\bQ})$ as in the previous subsection that also satisfies Assumption \ref{Assumption}.
\begin{defn2}Let $M\supset F_{\Pi}$ be a number field. A $\psi$-Whittaker function $W$ on $G$ is called to be \emph{defined over $M$} if it takes values in $M(\mu_{\infty})$ and satisfies
\[\sigma(W(g))=W(w(\kappa(\sigma))g),\]for all $g\in G(\bA_{f})$ and $\sigma\in\Gal(\ov{\bQ}/M)$, where $w(x)=\diag(x,1,x^{-1})$ and $\kappa:\Gal(\ov{\bQ}/\bQ)\rar\wh{\bZ}^{\times}$ is the cyclotomic character.
\end{defn2}Again, by Remark \ref{RKW}(1), $W$ is defined over $M$ if and only if $W(1)\in M^{\times}$.
\subsection{Whittaker periods via cohomological period integrals}
The purpose of this section is to prove the following
\begin{thm2}\label{U21W}Let $\varphi\in\Pi^{\new}=\Pi_{f}^{\new}\otimes\Pi_{\infty}^{\new}$. If $[\varphi]$ (see Definition \ref{CohclassDef}) is defined over a number field $F'\supset  F_{\Pi}$, then 
$\Lambda\left(1/2,BC(\Pi)\right)W_{\varphi}$ is a nontrivial $\psi$-Whittaker function defined over $F'$.
\end{thm2}
\begin{proof}
Let $C\in\bC^{\times}$ be such that $\frac{W_{\varphi}}{C}$ is defined over $F'$, which exists by Remark \ref{RKW}(2). We would like to show that
\[C\Lambda(1/2,BC(\Pi))\in F'{}^{\times}.\]First of all, this is indeed nonzero as $BC(\Pi)$ is cuspidal and tempered. 

We consider Gelbart--Piatetski-Shapiro integral representation of base change $L$-function. This is, roughly speaking, the period integral of $\varphi$ times an Eisenstein over an embedded modular curve. More precisely, there is an embedding of $H=\UU(1,1)$ 
into $G$. Let $B\subset H$ be the upper-triangular Borel. For $\Phi:\bA^{2}\rar\bC$ and $\chi$ a unitary Gr\"ossencharacter of the imaginary quadratic field $F$ used for the definition of unitary groups, we define an Eisenstein series on $H$ with respect to $B$, $E(h,\chi,\Phi,s)$. The same arithmeticity condition as Proposition \ref{Shim} applies, as the objects involved are elliptic modular forms.

Using the notation of Example \ref{ExampleReducible}(3), let $(m,n)=(a-b,b-c)$ be the Harish--Chandra character\footnote{Note that $G$ is now $\SU(2,1)$, so only the differences matter.} of $\Pi_{\infty}$. Let $\Gamma_{H}=H\cap\Gamma$, and let $i:X_{H}(\Gamma_{H})\hrar X_{G}(\Gamma)$ be the closed embedding of closed Shimura varieties, as before. Then, Gelbart--Piatetski-Shapiro's integral representation can be understood via a cup product pairing\[H^{1}(X_{G},V)\otimes_{\bC}H^{0}(X_{H},W)\xrar{(\id,i_{*})}H^{1}(X_{G},V)\otimes_{\bC}H^{1}(X_{G},W')\xrar{\cup} H^{2}(X_{G},V\otimes W')\xrar{S} H^{0}(X_{G},\cO)=\bC,\]where $\cup$ is the cohomological cup product, $S$ is a Serre duality pairing, normalized as in the proof of Theorem \ref{Sp4W} (namely, the pairing induced from a morphism of $K_{\infty}$-representations defined over $\ov{\bQ}$). Here, $W$ and $W'$ are certain automorphic vector bundles over $X_{H}$ and $X_{G}$, respectively, corresponding to $[V_{H}(|m-n|-1)]\otimes\omega_{H}(1,0)$ and $[\wt{L}_{1}]$, if we use the notation analogous to \cite[\S6]{LPSZ}. The reinterpretation, done in \cite{Oh}, of Gelbart--Piatetski-Shapiro's integral asserts that there is a coherent cohomology class $[E]\in H^{0}(X_{H}(\Gamma_{H}),W)$ defined over $FF_{\Pi_{f}}$, which corresponds to a nearly-holomorphic Eisenstein series $E$ under the Hodge splitting as in \cite[\S6.3]{LPSZ}, such that\[\langle[\varphi],[E]\rangle=C\Lambda\left(1/2,BC(\Pi)\right),\] where $BC$ means base-change. As $\langle,\rangle$ is defined over $\bQ$, the left hand side is in $F'$. On the other hand, the right hand side is nonzero as observed above. Thus, the desired statement follows.
\end{proof}
\subsection{Hodge-linear algebra}
We now conduct relevant calculations in Hodge-linear algebra to prove Theorem \ref{MainEvidence} for $G=\SU(2,1)$. For the sake of simplicity, we assume that we work with the case of Harish-Chandra character $(1,1,0)$; the same calculation yields the proof for general weights. Our goal is, as in \S\ref{Sp4Hodge}, to convert $\vol F^{1}H_{\dR}(\Ad M)$ into an expression that involves Deligne's periods. We will prove\begin{prop2}\label{SU21HodgeProp}For a motive $M$ associated with $\Pi$ in the sense of Conjecture \ref{ActualMotive}, we have

\[\vol F^{1}H_{\dR}(\Ad M)\sim_{\ov{\bQ}^{\times}}\frac{\pi^{3}c^{+}(BC(M))^{3/2}}{\delta(M_{\ov{\sigma}})^{1/2}}.\]\end{prop2}\begin{proof}In this case, the motive $M$  is a motive over $F$ with coefficients in $\bQ$. If we denote $\sigma:F\rar \bC$ by the preferred complex embedding, then as the minimal $K$-type is just $(1,1,1)$, by the recipe in \cite[2.3]{HLS},\[H_{B}(M_{\sigma})\otimes_{\bQ}\bC=\undertext{$H^{1,0}$}{\underbrace{\bC v_{1}\oplus\bC v_{2}}}\oplus\undertext{$H^{0,1}$}{\underbrace{\bC v_{3}}},\quad H_{B}(M_{\ov{\sigma}})\otimes_{\bQ}\bC=\undertext{$H^{-1,0}$}{\underbrace{\bC\ov{v}_{1}\oplus\bC\ov{v}_{2}}}\oplus\undertext{$H^{0,-1}$}{\underbrace{\bC\ov{v}_{3}}}.\]We can choose $v_{i}$ and $\ov{v}_{i}$ so that $F_{\infty}(v_{i})=\ov{v}_{i}$. Then, $BC(M):=\Res_{F/\bQ}M_{F}=M_{\sigma}\oplus M_{\ov{\sigma}}(-1)$ and $\Ad M=M_{\sigma}\otimes M_{\ov{\sigma}}$. The Deligne period $\delta$ of a motive, as before, is the determinant of the Betti-to-de Rham comparison map with respect to natural underlying $\bQ$-structures. Also, $c^{+}(BC(M))$ in this case would be the determinant of the map
\[H_{B}(BC(M))^{+}\otimes\bC\hrar H_{B}(BC(M))\otimes\bC\riso H_{\dR}(BC(M))\otimes\bC\thrar (H_{\dR}(BC(M))/F^{1}H_{\dR}(BC(M)))\otimes\bC,\]with respect to the natural underlying $\bQ$-structures.

Let $e_{1}^{+},e_{2}^{+},e_{3}^{+}$ be a $\bQ$-basis of $H_{B}(BC(M))^{+}$, $e_{1}^{-},e_{2}^{-},e_{3}^{-}$ be a $\bQ$-basis of $H_{B}(BC(M))^{-}$, $f_{1}f_{2},f_{3}\in F^{1}H_{\dR}(BC(M))$ be an $F$-basis, $g_{1},g_{2},g_{3}\in H_{\dR}(BC(M))/F^{1}H_{\dR}(BC(M))$ be an $F$-basis, and $\wt{g}_{1},\wt{g}_{2},\wt{g}_{3}\in H_{\dR}(BC(M))$ be lifts of $g_{1},g_{2},g_{3}$. We can further assume that $f_{1},f_{2}\in F^{1}H_{\dR}(M_{\sigma})$, $f_{3}\in F^{1}H_{\dR}(M_{\ov{\sigma}})$, $\wt{g}_{1},\wt{g}_{2}\in H_{\dR}(M_{\ov{\sigma}})$, $\wt{g}_{3}\in H_{\dR}(M_{\sigma})$. Given two $\bC$-bases of $H_{B}(M)\otimes\bC\cong H_{\dR}(M)\otimes\bC$, we can express the map into a matrix,
\[\begin{pmatrix} e_{1}^{+}&e_{2}^{+}&e_{3}^{+}&e_{1}^{-}&e_{2}^{-}&e_{3}^{-}\end{pmatrix}=\begin{pmatrix} f_{1}&f_{2}&f_{3}&\wt{g}_{1}&\wt{g}_{2}&\wt{g}_{3}\end{pmatrix}\begin{pmatrix} A&B\\C&D\end{pmatrix},\]for $A,B,C,D\in M_{3}(\bC)$. 

Under the canonical isomorphisms
\[H^{1,0}(BC(M))\cong F^{1}H_{\dR}(BC(M))\otimes\bC,\quad H^{0,1}(BC(M))\cong\frac{H_{\dR}(BC(M))}{F^{1}H_{\dR}(BC(M))}\otimes\bC,\]let \[f_{1,B},f_{2,B},f_{3,B}\in H^{1,0}(BC(M)),\qquad g_{1,B},g_{2,B},g_{3,B}\in H^{0,1}(BC(M)),\]be the images of $f_{1}, f_{2}, f_{3}, g_{1}, g_{2}, g_{3}$ under the above isomorphisms. Then, by the same argument as in \S\ref{Sp4Hodge},\[\begin{pmatrix} e_{1}^{+}&e_{2}^{+}&e_{3}^{+}&e_{1}^{-}&e_{2}^{-}&e_{3}^{-}\end{pmatrix}=\begin{pmatrix} f_{1,B}&f_{2,B}&f_{3,B}&{g}_{1,B}&{g}_{2,B}&{g}_{3,B}\end{pmatrix}\begin{pmatrix} A'&B'\\C&D\end{pmatrix}.\]As $c_{B}(f_{1,B}),c_{B}(f_{2,B}),c_{B}(f_{3,B})$ and $g_{1,B},g_{2,B},g_{3,B}$ are two $\bC$-bases of $H^{0,1}(BC(M))$, there is a system of linear relations
\[\begin{pmatrix}c_{B}(f_{1,B})&c_{B}(f_{2,B})&c_{B}(f_{3,B})\end{pmatrix}=\begin{pmatrix}g_{1,B}&g_{2,B}&g_{3,B}\end{pmatrix}X,\]for some $X\in\GL_{3}(\bC)$. Using the same technique as \cite[Lemma, \S8.2.1]{PV}, we see that $C=XA',D=XB'$, which similarly implies that \[c^{+}(BC(M))^{2}\sim_{\ov{\bQ}^{\times}}\delta(BC(M))\det X,\]where we used $c^{+}(BC(M))\sim_{\ov{\bQ}}c^{-}(BC(M))$ due to the decomposition $BC(M)=M_{\sigma}\oplus M_{\ov{\sigma}}(-1)$ as in \cite[(8)]{HarrisLin}. Furthermore, as $g_{1,B},g_{2,B}\in H^{0,1}(M_{\sigma})$ and $g_{3,B}\in H^{0,1}(M_{\ov{\sigma}}(-1))$, it turns out that $X$ is of form \[\begin{pmat}Y&\\&\mu\end{pmat},\quad Y\in\GL_{2}(\bC),\mu\ne0.\]
We know that $(\vol F^{1}H_{\dR}(\Ad M))^{2}\sim_{\ov{\bQ}^{\times}}\lambda$, where $\varphi(c_{B}(v^{+}))=\lambda v^{-}$, $v^{+}=(f_{1}\otimes f_{3})\wedge (f_{2}\otimes f_{3})$ and $v^{-}=(g_{1}\otimes g_{3})\wedge (g_{2}\otimes g_{3})$. This implies that $\lambda=\mu^{2}\det Y=\mu\det X$. 

Now, as in \cite[\S1.2]{HarrisPeriod2}, consider the determinant motive $\det(M)$. An $F$-rational basis vector of $H_{\dR}(\det(M_{\sigma}))$ can be taken as $v_{\sigma}:=f_{1}\wedge f_{2}\wedge \wt{g}_{3}$, and similarly $v_{\ov{\sigma}}:=\wt{g}_{1}\wedge\wt{g}_{2}\wedge f_{3}$ can be taken as an $F$-rational basis vector of $H_{\dR}(\det(M_{\ov{\sigma}}))$. On the other hand, if we take $e_{\sigma}$ to be a $\bQ$-rational basis vector of $H_{B}(\det(M_{\sigma}))$, then $F_{\infty}(e_{\sigma})=:e_{\ov{\sigma}}$ is a $\bQ$-rational basis vector of $H_{B}(\det(M_{\ov{\sigma}}))$. Then $e_{\sigma}=\delta(M_{\sigma})v_{\sigma}$ and $e_{\ov{\sigma}}=\delta(M_{\ov{\sigma}})v_{\ov{\sigma}}$, so
\[c_{B}(v_{\sigma})=\delta(M_{\sigma})^{-1}e_{\ov{\sigma}}.\]On the other hand,
\[c_{B}(v_{\sigma})=F_{\infty}(v_{\sigma})=(\det Y\cdot \mu^{-1}) v_{\ov{\sigma}},\]so 
\[\delta(M_{\ov{\sigma}})=\mu^{-1}\delta(M_{\sigma})\det Y.\]
On the other hand, due to the polarization, we have \cite[(1.2.5)]{HarrisPeriod2}, 
\[\delta(M_{\ov{\sigma}})=\delta(M_{\sigma})^{-1}(2\pi i)^{-6}.\]Thus,
\[\pi^{6}\delta(M_{\sigma})^{2}\det Y\sim_{\ov{\bQ}^{\times}}\mu,\]so 
\[\mu^{2}\sim_{\ov{\bQ}^{\times}}\pi^{6}\delta(M_{\sigma})^{2}\det X\sim_{\ov{\bQ}^{\times}}c^{+}(BC(M))^{2}\pi^{6}\frac{\delta(M_{\sigma})}{\delta(M_{\ov{\sigma}})}\sim_{\ov{\bQ}^{\times}}c^{+}(BC(M))^{2}\pi^{12}\delta(M_{\sigma})^{2},\]or
\[\mu\sim_{\ov{\bQ}^{\times}}\pi^{6}\delta(M_{\sigma})c^{+}(BC(M)),\]which implies that
\[(\vol F^{1}H_{\dR}(\Ad M))^{2}\sim_{\ov{\bQ}^{\times}}\frac{\pi^{6}}{\delta(M_{\ov{\sigma}})}c^{+}(BC(M))^{3},\]as desired.
\end{proof}
\subsection{Completion of the proof}
\begin{proof}[Proof of Theorem \ref{MainEvidence} for $\SU(2,1)$]
Let $f_{\hol},f_{\gen}$ be newforms (see Assumption \ref{Assumption}) in $\Pi_{f}\otimes\Pi_{\hol}$ and $\Pi_{f}\otimes\Pi_{\gen}$, respectively, such that $[f_{\hol}],[f_{\gen}]$ are defined over $\ov{\bQ}$. Then, under the assumptions of (\ref{SU21A}) of Assumption \ref{Assumption}\footnote{There is an extra condition on large residue characterstic in \emph{loc. cit.}, but this restriction is recently removed by \cite[Theorem 1]{BP}.}, \cite[Theorem 1.2]{Wei}\footnote{This is a special case of the refined Gan--Gross--Prasad conjecture for Fourier--Jacobi periods (see e.g. \cite[\S1.1]{Xue}), which is also often referred as the Ichino--Ikeda conjecture.} implies that
\[\tag{F}\frac{\langle f_{\hol},f_{\hol}\rangle_{P}}{|FJ(f_{\hol})|^{2}}\sim_{\ov{\bQ}^{\times}}\pi^{-2}\frac{L(1,\Pi,\Ad)}{L(1/2,BC(\Pi))},\]where $FJ(f_{\hol})$ is the special Fourier--Jacobi period of $f_{\hol}$, defined by
\[FJ(f_{\hol})=\int_{\SU(1,1)(\bQ)\bs\SU(1,1)(\bA)}f_{\hol}(h)dh,\]integrated against the Tamagawa measure. This period is in turn expressed as an inner product of (algebraic) theta functions,
\[FJ(f_{\hol})\sim_{\ov{\bQ}^{\times}}\langle a(0,f_{\hol})(v),1\rangle,\]where $a(0,f_{\hol})$ is the zero-th Fourier--Jacobi coefficient of $f_{\hol}$ and $1$ is the constant function (regarded as a trivial theta function). Note that the nonvanishing of $FJ(f_{\hol})$ (and thus $a(0,f_{\hol})$) is also a part of the content of \cite[\S1.1]{Wei}.

By \cite{LanFJ}, $a(0,f_{\hol})$ is identified with the algebraic Fourier--Jacobi coefficient, and in particular is in $\ov{\bQ}^{\times}$, as $[f_{\hol}]$ is defined over $\ov{\bQ}$. Thus, we have $\langle f_{\hol},f_{\hol}\rangle_{P}\sim_{\ov{\bQ}^{\times}}\pi^{-2}\frac{L(1,\Pi,\Ad)}{L(1/2,BC(\Pi))}$.

On the other hand, we exploit the fact that $H^{1}$ is the middle degree of the Shimura variety. Note that the infinity type of $\ov{f}_{\gen}$ is a generic LDS as well (with different infinitesimal character from $f_{\gen}$). Suppose $C\in\bC^{\times}$ is a constant where $C[\ov{f}_{\gen}]$ is defined over $\ov{\bQ}$. Then,
\begin{align*}\tag{G}
\langle f_{\gen},f_{\gen}\rangle_{P}&=(2\pi i)^{2}\langle [f_{\gen}],[\ov{f}_{\gen}]\rangle_{\coh}\\&\sim_{\ov{\bQ}^{\times}}\pi^{2}C^{-1},
\end{align*}where the cohomological cup product pairing $\langle-,-\rangle_{\coh}$ is the Serre duality pairing induced from a $\ov{\bQ}$-morphism of algebraic $K_{\infty}$-representations $V\otimes\Hom(V,\fg^{-1,1})\rar\fg^{-1,1}$.

By Theorem \ref{U21W}, \[\Lambda(1/2,BC(\Pi))W_{C\ov{f}_{\gen}}=C\Lambda(1/2,BC(\Pi))W_{\ov{f}_{\gen}},\]is defined over $\ov{\bQ}$.
We now invoke the Lapid--Mao conjecture (Conjecture \ref{LMConj}), which says
\[\langle f_{\gen},f_{\gen}\rangle_{P}\sim_{\bQ^{\times}}
|W_{f_{\gen}}(1)|^{2}\frac{L(1,\Pi,\Ad)}{|W_{\infty}(1)|^{2}}.\]By Remark \ref{RKW}(1), $W_{\ov{f}_{\gen}}(1)\ne0$, and $C\Lambda(1/2,BC(\Pi))W_{\ov{f}_{\gen}}(1)\in\ov{\bQ}^{\times}$ as observed in the beginning of the section.
Thus, we have\[|W_{f_{\gen}}(1)|^{2}\frac{L(1,\Pi,\Ad)}{|W_{\infty}(1)|^{2}}\overtext{Lapid--Mao}{\sim_{\ov{\bQ}^{\times}}}\langle f_{\gen},f_{\gen}\rangle_{P}\overtext{(F)\quad}{\sim_{\ov{\bQ}^{\times}}}\pi^{2}\Lambda(1/2,BC(\Pi))\ov{W_{f_{\gen}}(1)}.\]By using $W_{\ov{f}_{\gen}}(1)=\ov{W_{f_{\gen}}(1)}$, we have
\[\tag{H}W_{f_{\gen}}(1)\sim_{\ov{\bQ}^{\times}}\pi^{2}\frac{|W_{\infty}(1)|^{2}\Lambda(1/2,BC(\Pi))}{L(1,\Pi,\Ad)}.\]On the other hand, Theorem \ref{U21W} says
\[W_{f_{\gen}}(1)\sim_{\ov{\bQ}^{\times}}\frac{1}{\Lambda(1/2,BC(\Pi))},\]which gives an extra relationship that we can utilize; namely, \[\Lambda(1,\Pi,\Ad)\sim_{\ov{\bQ}^{\times}}\Lambda(1/2,BC(\Pi))^{2}|W_{\infty}(1)|^{2}\pi^{2}.\] For example, by combining (F) and (G), we have
\[\tag{I}\langle f_{\gen},f_{\gen}\rangle_{P}\sim_{\ov{\bQ}^{\times}}\pi^{4}\frac{|W_{\infty}(1)|^{2}\Lambda(1/2,BC(\Pi))^{2}}{L(1,\Pi,\Ad)}.\]\begin{rmk2}Note that this is already observed in \cite[Corollary 1.3.5]{HarrisPeriod2}, under the assumption of Deligne's conjectures. In particular, assuming Deligne's conjectures, we deduce that $|W_{\infty}(1)|\sim_{\ov{\bQ}^{\times}\pi^{\bZ}}1$.\end{rmk2}

Combining these, we get
\begin{align*}\frac{\langle f_{\hol},f_{\hol}\rangle_{P}}{\langle f_{\gen},f_{\gen}\rangle_{P}}&\overtext{Lapid--Mao + (H)}{\sim_{\ov{\bQ}^{\times}}}\frac{\pi^{-2}\frac{L(1,\Pi,\Ad)}{L(1/2,BC(\Pi))}}{\pi^{4}\frac{|W_{\infty}(1)|^{2}\Lambda(1/2,BC(\Pi))^{2}}{L(1,\Pi,\Ad)}}\\&\sim_{\ov{\bQ}^{\times}}\frac{1}{|W_{\infty}(1)|^{2}}\frac{L(1,\Pi,\Ad)^{2}}{L(1/2,BC(\Pi))^{3}}\\&\sim_{\ov{\bQ}^{\times}}\frac{\pi^{18}}{|W_{\infty}(1)|^{2}}\left(\frac{L_{\infty}(1,\Pi,\Ad)}{L_{\infty}(0,\Pi,\Ad)}\cdot\frac{L(1,\Pi,\Ad)}{L(1/2,BC(\Pi))^{3/2}}\right)^{2},\\&\sim_{\ov{\bQ}^{\times}}\frac{\pi^{21}}{|W_{\infty}(1)|^{2}}\left(\frac{L_{\infty}(1,\Pi,\Ad)}{L_{\infty}(0,\Pi,\Ad)}\cdot\frac{L(1,\Pi,\Ad)}{\vol F^{1}H_{\dR}(\Ad M)}\right)^{2},\end{align*}by Deligne's conjecture applied to Proposition \ref{SU21HodgeProp}. We are done, as the paranthesized term is the RHS of Conjecture \ref{Period}.\end{proof}
\begin{rmk2}\label{SU21Whittaker}
Similarly to Remark \ref{Sp4Whittaker}, $W_{\infty}(1)$ is given by the inverse Mellin transform of the formulae given in \cite[Theorem 5.5]{KO}, which the author at the moment is unable to compute.
\end{rmk2}
\section{Towards motivic action conjecture for rationality of classes}\label{MotivicRationalSect}
The original motivic action conjecture of \cite{PV} involves the rational structure of singular cohomology. To derive a similar conjecture, we would have to come up with a way to normalize all the (choice of newforms of) automorphic representations at once. Recall that, in the case of modular forms, this is done by using complex conjugation. Unfortunately, so far there is no general construction of an operation that can move between different infinity types. We will tentatively name such an operation a \emph{generalized complex conjugation}. Approaching the generalized complex conjugations using the theory of cycle spaces and Penrose transform will be the subject of the author's forthcoming work. For now, we will have to content ourselves with a preliminary analysis on what a generalized complex conjugation should be. 

On the other hand, as the name suggests, the usual complex conjugation can be used when the symmetric space is the upper half plane. More generally, if the symmetric space is a product of upper half planes (e.g. in the case of Hilbert modular varieties), then the complex conjugations with respect to each variable would be a good candidate for generalized complex conjugation; these are called \emph{partial complex conjugations} \cite{HarrisHilb} in the literature. In this case, we can deduce a precise conjecture on the $\ov{\bQ}$-rational structure of coherent cohomology, from the generalities of Appendix \ref{ExtSection}. There is an existing work of Horawa exactly on this problem \cite{Horawa}, and we will compare our conjecture with that of \emph{op. cit.} In particular, we observe that the numerical evidences given in \emph{op. cit.} are compatible with both conjectures. 
\subsection{The case of Hilbert modular forms: comparison with \cite{Horawa}}\label{HorawaSect}

The work \cite{Horawa} states a similar conjecture, Conjecture 3.21 of \emph{op. cit.}, on what archimedean motivic action should be for Hilbert modular forms of partial weight one. It uses the partial complex conjugation, which utilizes the fact that every (limit of) discrete series for $\SL_{2}(\bR)^{d}$ is holomorphic or antiholomorphic in each variable. 
\begin{defn2}[(Partial complex conjugations)]Let $F$ be a totally real field of degree $d$ and $\varphi$ be a holomorphic automorphic form for $G=\Res_{F/\bQ}\GL_{2,F}$, seen as a holomorphic function on the symmetric space for $G(\bR)$, $(\bC-\bR)^{d}$, of weight $(k_{1},\cdots,k_{d};r)$ where $k_{i}\ge1$, $k_{i}\equiv r(\MOD 2)$ for $i=1,\cdots,d$. For $I\subset\lbrace1,\cdots,d\rbrace$, $\varphi^{I}$ is the automorphic form for $G=\Res_{F/\bQ}\SL_{2,F}$, defined by $\varphi^{I}(g)=\varphi(gJ^{I})$ for $J^{I}=(J_{1}^{I},\cdots,J_{d}^{I})$ given by $J_{j}^{I}=\begin{psmat}-1&0\\0&1\end{psmat}$ if $j\in I$ and $J_{j}^{I}=\id_{2}$ if $j\notin I$. This is called the \emph{partial complex conjugation}.
\end{defn2}The main 
conjecture of op. cit., \cite[Conjecture 3.21]{Horawa}, describes the rationality of cohomology classes in terms of partial complex conjugations. In particular, it implies that the decomposition $H^{i}(X)[\Pi_{f}]=H^{i}(\fp,K;\omega\otimes I(\lambda))\otimes\Pi_{f}^{\Gamma}=\left(\bigoplus_{\pi\in\fP_{\lambda}}H^{i}(\fp,K;\omega\otimes\pi)\right)\otimes\Pi_{f}^{\Gamma}$ descends to $\ov{\bQ}$. 

We have not been successful in approaching the conjecture.
On the other hand, based on the materials developed in Appendix \ref{ExtSection}, we suggest a slightly different conjecture. We use the language of \cite[\S3]{Horawa} belwo.
\begin{conj2}[(Motivic action conjecture for Hilbert modular varieties; compare with {\cite[Conjecture 3.21]{Horawa}})]\label{OurHilb}Let $f$ be a parallel weight one form. For each $u\in U_{f}^{\vee}$, let $u_{i}\in(\Ad^{0}M\otimes_{\iota}\bC)^{\sigma_{i}c_{0}\sigma_{i}^{-1}}\cong\bC$ be the $\sigma_{i}$-component of $U_{f}^{\vee}\otimes_{\iota}\bC$ as in \cite[Proposition 3.2]{Horawa}, where the isomorphism is given by the natural $\ov{\bQ}$-structure on $\fsl_{2}^{d}$. Then, for every $u\in U_{f}^{\vee}$ not in the kernel of the pairing of \cite[Lemma 3.1]{Horawa},\[2\pi i\sum_{i=1}^{d}\frac{\omega_{f}^{\lbrace i\rbrace}}{\log(|\tau\otimes\iota(u)|)}\in H^{1}(X(\Gamma),\omega),\]defines a cohomology class over $\ov{\bQ}$.
\end{conj2}
Conjecture \ref{OurHilb} suggests that  the decomposition of coherent cohomology as a $(\fp,K)$-cohomology of archimedean representations is not in general $\ov{\bQ}$-rational. This is compatible with the corrected version of Conjecture 3.21 of \emph{op. cit.}.

In \cite[\S5]{Horawa}, a numerical evidence in favor of the conjecture of \cite{Horawa} is given for base change forms in the case of Hilbert modular forms for real quadratic fields. We claim that, for such Hilbert modular forms, the two conjectures coincide:
\begin{prop2}
Let $f$ be a Hilbert modular eigen-cuspform of parallel weight one for a real quadratic field $F$. If $f$ is a base change form, then \cite[Conjecture 3.21]{Horawa} implies Conjecture \ref{OurHilb}.
\end{prop2}
\begin{proof}Let $\sigma_{1},\sigma_{2}$ be the two real embeddings of $F$, and suppose $f$ is a base change form of $f_{0}$. Indeed, it is shown in \cite[Corollary 5.2]{Horawa} that the space of $\ov{\bQ}$-Stark units for $f$, $U_{f}\otimes\ov{\bQ}$, is naturally isomorphic to $(U_{f_{0}}\otimes\ov{\bQ})^{\oplus 2}$, and the decomposition is compatible with the Beilinson regulator. In particular, the $\ov{\bQ}$-vector space spanned by the log of Stark units of ${f}$ is exactly the $\ov{\bQ}$-vector space spanned by the log of Stark units of $f_{0}$. In particular, both Conjecture \ref{OurHilb} and \cite[Conjecture 3.21]{Horawa} are equivalent to the statement that \[\ov{\bQ}\omega_{f}^{\lbrace1\rbrace}\oplus\ov{\bQ}\omega_{f}^{\lbrace2\rbrace}=\frac{\log(u_{f_{0}})}{2\pi i}H^{1}(X(\Gamma)_{\ov{\bQ}},\omega),\] where $u_{f_{0}}\in U_{f_{0}}$.
\end{proof}
Indeed, a base change form satisfies an extra symmetry with respect to the ``change of two upper half planes'', namely \[\bH^{2}\xrar{(x,y)\mapsto(y,x)}\bH^{2},\]and this extra symmetry guarantees that the $\ov{\bQ}$-splitting of our form is compatible with the $\ov{\bQ}$-splitting of the form in \cite[Conjecture 3.21]{Horawa}.

\subsection{Desiderata for generalized complex conjugations}\label{GCCSection}
To have a normalized choice of newforms simultaneously, we would like a certain way to relate different newforms. Example \ref{SL2} suggests that the complex conjugation should play a role in the tentative statement of the full motivic action conjecture regarding rationality of cohomology classes. On the other hand, the complex conjugation can go back and forth between only two types of LDS's. For example, it sends a vector in the ``holomorphic LDS'' to that in the ``anti-holomorphic LDS.'' Since a general motivic action conjecture involves many more LDS's, we suggest that there are \emph{generalized complex conjugations} that can go between any of $\pi\in\fP_{\lambda}$. 

We denote a generalized complex conjugation, sending an automorphic form $v\in\Pi\otimes A_{C}(\lambda)$ to another automorphic form in $\Pi\otimes A_{C'}(\lambda)$, by $c_{C,C'}$. There are several desired properties:
\begin{itemize}
\item $c_{C,C'}$ is $\bC$-linear,
\item If $C_{\hol},C_{\mathrm{antihol}}$ are in $\fC_{\lambda}$, $c_{C_{\hol},C_{\mathrm{antihol}}}(f)=\ov{f}$,
\item $\langle f,f\rangle_{P}=\langle c_{C,C'}(f),c_{C,C'}(f)\rangle_{P}$ (Condition (B)),
\item $c_{C',C''}\circ c_{C,C'}=c_{C,C''}$ and $c_{C,C}=\id$,
\item $c_{C_{1}'\times C_{2}',C_{1}\times C_{2}}=(c_{C_{1}',C_{1}},c_{C_{2}',C_{2}})$, for $G(\bR)=G_{1}(\bR)\times G_{2}(\bR)$.
\end{itemize}
It is still unclear how to formulate a set of conditions which will uniquely characterize $c_{C,C'}$'s. Although the nature of generalized complex conjugations still remains mysterious, using the ideas of Penrose transform and its related geometry, it could be possible to construct the purported generalized complex conjugations, following the suggestion by Joseph Wolf. This is the subject of the author's forthcoming work. 
\begin{rmk2}
The last bullet point suggests that the \emph{partial complex conjugation} (e.g. \cite{HarrisHilb}, \cite{ShimuraHilb}) serves the role of generalized complex conjugations in the case of Hilbert modular forms
. Unfortunately, for $C_{\hol}$ and $C_{\mathrm{anithol}}$ to be both in $\fC_{\lambda}$, $\lambda$ has to be orthogonal to all compact roots, and this is allowed only if there is no compact root (as we exclude degenerate limit of discrete series from our discussion). Thus, we cannot use the usual complex conjugation besides when the associated symmetric space is a product of several upper half planes.
\end{rmk2}

\begin{rmk2}
It is a relatively well-accepted technique in the case of unitary groups to use theta correspondence to move between different infinity types, as suggested by the recipe of \cite{Prasad}. Indeed, for a unitary group, \cite{HLS} proves that the theta correspondences and character twists act transitively upon the full Vogan $L$-packet (see \cite{VoganLLC}). However, due to the idiosyncrasies of the recipe for the theta correspondence, it is still unclear whether the theta correspondence should be the generalized complex conjugation in this case.
\end{rmk2}
\begin{rmk2}
We also speculate that this is the archimedean version of \emph{excursion operators} (first appeared in the work of \cite{Lafforgue} on the global Langlands correspondence over function fields, and extended to the context of mixed characteristic local Langlands via Kottwitz's conjecture, e.g. \cite{FarguesMantovan}, \cite{RapoportViehmann}). Indeed, the isotypic decomposition with respect to the excursion algebra canonically decomposes the automorphic spectrum into $L$-indistinguishable pieces, which would mean that the excursion operators can go around different members of an $L$-packet. 
\end{rmk2}

Under the hypothesis on existence of generalized complex conjugations $c_{C,C'}$, we can formulate the motivic action conjecture in the Shimura variety context in its full form.
\begin{conj2}\label{Main}
Let $\lambda,\Pi$ as in Conjecture \ref{Metric}. Let $f_{h}\in\Pi_{f}^{\new}\otimes\pi_{h}^{\new}$ be a newform such that $[f_{h}]$ (see Definition \ref{CohclassDef}) is defined over $\ov{\bQ}$. Let $\cE_{\lambda}\cong\Ext^{1}_{(\fp,K)}(I(\lambda),\pi_{h})$ be defined such that $\bC\alpha_{i}\cong\Ext^{1}_{(\fp,K)}(\pi_{\lbrace 1,\cdots,n_{\lambda}\rbrace-\lbrace i\rbrace},\pi_{h})$ sends $1\alpha_{i}$ to the homomorphism $c_{C_{h},C_{\lbrace 1,\cdots,n_{\lambda}\rbrace-\lbrace i\rbrace}}(f_{h})\mapsto f_{h}$ (using the identification from Proposition \ref{ExtHom}). Then, for $v\in H_{M}^{1}((\Ad^{*}\Pi)_{\cO_{E}},\ov{\bQ}(1))\subset H_{\sD}^{1}((\Ad^{*}\Pi)_{\bR},\bC(1))$, $a(v)\cdot f_{l}$ defines a coherent cohomology class $[a(v)\cdot f_{l}]\in H^{i_{l}+1}(X_{G}(\Gamma),[V_{A(\lambda)}])$ that is defined over $\ov{\bQ}$.
\end{conj2}One can easily state a similar conjecture for the action of $\bigwedge^{*}H_{M}^{1}((\Ad\Pi)_{\cO_{E}},\ov{\bQ}(1))$, but it is no deeper than the conjecture stated above.

There is basically one known case of what generalized complex conjugation should be, and it is the case of Hilbert modular forms.
\appendix
\section{Beilinson's conjecture over a general number field}\label{BCSection}
In this section, we recall the statement of \emph{Beilinson's conjecture} we will need in the paper. A usual formulation of the conjecture involves motives over $\bQ$ with coefficients in $\bQ$, but we would have to relax both to be arbitrary number fields. A standard reference of this matter is \cite[\S6]{Ramakrishnan}.
\subsection{Chow motives}We recall the definition of Chow motives over a number field $k$, $\sM_{k,\rat}$ in \cite[\S2.1.1]{PV}. defined by cohomological correspondences up to rational equivalence. If $k$ is a number field, then for a Chow motive $M\in\sM_{k,\rat}$, there are the following cohomology theories, motivated by the cohomology theories of smooth proper $k$-varieties.
\begin{itemize}
\item For each prime $\ell$, there is $\ell$-adic cohomology $H^{i}(M_{\ov{K}},\bQ_{\ell}(r))$, which is a finite-dimensional $\ell$-adic representation of $\Gal(\ov{K}/K)$.
\item For each embedding $\sigma:k\hrar\bC$, there is Betti cohomology $H^{i}_{B}(M_{\sigma},\bQ(r))$, which is a pure $\bQ$-Hodge structure. If $\sigma$ is a real embedding, it is equipped with the \emph{infinite Frobenius} $\Fr_{\sigma,\infty}$. On $H^{i}_{B}(M_{\sigma},\bQ(r))\otimes_{\bQ}\bC$, the involution $\Fr_{\sigma,\infty}\otimes c_{B}$ preserves the Hodge decomposition, where $c_{B}$ is the complex conjugation on the second factor. If $\sigma$ is a complex embedding, $\Fr_{\sigma,\infty}$ is rather an isomorphism of $\bQ$-vector spaces
\[\Fr_{\sigma,\infty}:H_{B}^{i}(M_{\sigma},\bQ(r))\riso H^{i}_{B}(M_{\ov{\sigma}},\bQ(r)).\]
\item There is de Rham cohomology $H^{i}_{\dR}(M)(r)$, which is a finite-dimensional $k$-vector space, equipped with a decreasing filtration $F^{k}H_{\dR}^{i}(M)(r)$.
\end{itemize}
We will not care much about $\ell$-adic realization, as it plays no role in the paper. For each embedding $\sigma:k\hrar \bC$, there is a comparison isomorphism
\[\comp_{\sigma}:H_{B}^{i}(M_{\sigma},\bQ(r))\otimes_{\bQ}\bC\cong H_{\dR}^{i}(M)(j)\otimes_{k,\sigma}\bC,\]such that $\oplus_{p\ge k}H^{p,q}$ corresponds to $(F^{k}H_{\dR})\otimes\bC$. If $\sigma$ is a real embedding, $\Fr_{\sigma,\infty}\otimes c_{B}$ corresponds to $1\otimes c_{\dR}$, where $c_{\dR}$ is the complex conjugation on the second factor of $H_{\dR}^{i}(M)(j)\otimes_{k,\sigma}\bC$. If $\sigma$ is a complex embedding, $\Fr_{\sigma,\infty}\otimes c_{B}$ corresponds to $1\otimes c_{\dR}$ in the sense that there is a commutative diagram
\[\xymatrix{H_{B}^{i}(M_{\sigma},\bQ(r))\otimes_{\bQ}\bC\ar[r]^-{\comp_{\sigma}}\ar[d]_-{\Fr_{\sigma,\infty}\otimes c_{B}}&H_{\dR}^{i}(M)(j)\otimes_{k,\sigma}\bC\ar[d]^-{1\otimes c_{\dR}}\\H_{B}^{i}(M_{\ov{\sigma}},\bQ(r))\otimes_{\bQ}\bC\ar[r]_-{\comp_{\ov{\sigma}}}&H_{\dR}^{i}(M)(j)\otimes_{k,\ov{\sigma}}\bC}\]
Another key player is the Deligne cohomology. For a complex smooth projective variety $X$ and a subring $A\subset\bC$ invariant under complex conjugation, the Deligne cohomology $H_{\sD}^{i}(X,A(r))$ is defined as the hypercohomology of the complex
\[A(r)_{\sD}:A(r)\rar\cO_{X}\xrar{d}\Omega_{X}^{1}\xrar{d}\cdots\xrar{d}\Omega_{X}^{r-1},\]regarded as a complex of analytic sheaves. This admits a complex conjugation of coefficients, denoted $c_{\sD}$, which is induced from the complex conjugation on $A(r)_{\sD}$. Similarly, there is infinite Frobenius for Deligne cohomology.

In the formulation of Beilinson's conjecture over $\bQ$, a central role is played by the cohomology of $M_{\bR}$. It plays the same role in Beilinson's conjecture over general number fields, even though it may sound peculiar to consider the base-change of $M$ to $\bR$ even if $k$ is not a real field. 
\begin{defn2}
Given a subfield $A\subset\bC$ stable under complex conjugation, define
\[H_{\dR}^{i}(M_{\bR})(r)=\left(\bigoplus_{\sigma:k\hrar\bC}H_{\dR}^{i}(M)(r)\otimes_{k,\sigma}\bC\right)^{1\otimes c_{\dR}}\]
\[H_{B}^{i}(M_{\bR},A(r))=\left(\bigoplus_{\sigma:k\hrar\bC}H_{B}^{i}(M_{\sigma},\bQ(r))\otimes_{\bQ}A\right)^{\oplus_{\sigma}\Fr_{\sigma,\infty}\otimes c_{B}}\]
\[H_{\sD}^{i}(M_{\bR},A(r))=\left(\bigoplus_{\sigma:k\hrar\bC}H_{\sD}^{i}(M_{\sigma},A(r))\right)^{\oplus_{\sigma}\Fr_{\sigma,\infty}\otimes c_{\sD}}.\]
\end{defn2}
Concretely, cohomology of $M_{\bR}$ is the part fixed by $\Fr_{\infty}\otimes c_{B}=1\otimes c_{\dR}$ via Betti-de Rham comparison isomorphism. Furthermore, there is a Chern class map\[r_{\sD}:H_{M}^{i}(M,\bQ(r))\rar H_{\sD}^{i}(M_{\bR},\bR(r)).\]
The source of Chern class map is too large, and we choose a subspace $H_{M}^{i}(M_{\cO_{k}},\bQ(r))\subset H_{M}^{i}(M,\bQ(r))$ consisting of classes that ``extend to a good proper model of $M$ over $\cO_{k}$''. If $M=h(X)$ for a smooth proper $k$-variety $X$, which has a regular proper model $\fX$ over $\cO_{k}$, then
\[H_{M}^{i}(M_{\cO_{k}},\bQ(r))=H_{M}^{i}(M,\bQ(r))\cap(\im(K_{2r-i}\fX\rar K_{2r-i}X)\otimes\bQ),\]where the latter image of $K$-theory groups is the image via the Chern class characters, and this definition is independent of choice of $\fX$. In general, using alterations, Scholl defined this subspace in \cite[Theorem 1.1.6]{Scholl} and showed that this is a unique way to assign subspaces satisfying various natural properties. The restriction of Chern class charcater into the integral subspace,
\[r_{\sD}:H^{i}_{M}(M_{\cO_{k}},\bQ(r))\rar H_{\sD}^{i}(M_{\bR},\bR(r)),\]is called the Beilinson regulator.
\subsection{Beilinson's conjecture for Chow and Grothendieck motives}
From now on, we assume that $r\ge\frac{i+1}{2}$\footnote{This is equivalent to that the weight of $M$, $i-2r$, is negative. This is not a restriction as one can always reduce to this case possibly after using functional equation.}. Beilinson's conjecture is formulated using fundamental exact sequences, which we review. From the definition of Deligne cohomology, for a complex smooth projective variety $X$, there is a long exact sequence
\[\cdots\rar H^{i-1}_{B}(X,A(r))\rar \frac{H^{i-1}_{B}(X,\bC)}{F^{r}H_{\dR}^{i-1}(X)}\rar H_{\sD}^{i}(X,A(r))\rar H_{B}^{i}(X,A(r))\rar\cdots.\]
As this long exact sequence intertwines the involutions that define the cohomology of $M_{\bR}$, for a Chow motive $M$ defined over a number field $k$, 
\[H_{\sD}^{i}(M_{\bR},\bR(r))=\frac{H_{B}^{i-1}(M_{\bR},\bC)}{F^{r}H_{\dR}^{i-1}(M_{\bR})+H_{B}^{i-1}(M_{\bR},\bR(r))}=\frac{H_{B}^{i-1}(M_{\bR},\bR(r-1))}{F^{r}H_{\dR}^{i-1}(M_{\bR})},\]which gives rise to two fundamental exact sequences,
\[0\rar F^{r}H_{\dR}^{i-1}(M_{\bR})\rar H_{B}^{i-1}(M_{\bR},\bR(r-1))\rar H_{\sD}^{i}(M_{\bR},\bR(r))\rar0,\]
\[0\rar H_{B}^{i-1}(M_{\bR},\bR(r))\rar \frac{H_{\dR}^{i-1}(M_{\bR})}{F^{r}H_{\dR}^{i-1}(M_{\bR})}\rar H_{\sD}^{i}(M_{\bR},\bR(r))\rar0.\]
The first two entries of the two fundamental exact sequences as above have natural $\bQ$-structures, yielding the $\bQ$-structure on $\det H_{\sD}^{i}(M_{\bR},\bR(r))$. Let $\sR$ be the one from the first sequence, and $\sD\sR$ be the one from the second exact sequence.
\begin{conj2}[(Beilinson's conjecture)]
Suppose either $r>\frac{i}{2}+1$, or $r=\frac{i}{2}+1$ and there is no Tate cycle, namely $H_{\et}^{i}(M_{\ov{k}},\bQ_{\ell}(i/2))^{\Gal(\ov{k}/k)}=0$. Then, the following hold.
\begin{enumerate}
\item The Beilinson regulator
\[r_{\sD}:H_{M}^{i+1}(M_{\cO_{k}},\bQ(r))\otimes\bR\rar H_{\sD}^{i+1}(M_{\bR},\bR(r)),\]is an isomorphism.
\item Let $\sM$ be the $\bQ$-structure defined on $\det H_{\sD}^{i+1}(M_{\bR},\bR(r))$ via the Beilinson regulator and the $\bQ$-structure of the motivic cohomology. Then, 
\[\sM=L(h^{-i}(M^{\vee}),1-r)^{*}\sR=L(h^{i}(M),r)\sD\sR.\]
\end{enumerate}
\end{conj2}\begin{rmk2}
The motives we will be working with come from automorphic representations in the sense of \cite[\S4.3.3]{Clozel}, and for that purpose, we would rather like to work with \emph{Grothendieck motives}, where the equivalence relation used is numerical equivalence, which is equivalent to homological equivalence under the Standard Conjecture D, which we have to assume. Fortunately, \cite[\S2.1.9]{PV} works with arbitrary base field, so a similar set of assumptions would naturally lead to Beilinson's conjecture for Grothendieck motives.
\end{rmk2}
\section{Deligne cohomology and Lie algebra cohomology}\label{ExtSection}
In this appendix, we develop a representation theoretic background parallel to \cite[\S2-\S4]{PV}. This is to correctly guess the motivic action conjectures, namely Conjecture \ref{Metric} and Conjecture \ref{Main}.

\subsection{Nondegenerate limit of discrete series as constituents of reducible principal series}\label{A1}
We need to understand what kinds of infinity types can appear in the coherent cohomology of Shimura varieties, given the finite part and the coefficient. Fortunately, the $(\fp,K)$-cohomology of unitary $G(\bR)$-representations is computed by Vogan--Zuckerman. As far as $G(\bR)$-representations are concerned, we will be only interested in discrete series, or in general nondegenerate limits of discrete series.

For a (L)DS, its $L$-packet consists of all (L)DS with the same infinitesimal character. Thus each such $L$-packet is consisted of $|W_{G}|/|W_{K}|$ elements, and in particular, upon choosing a system of positive roots $\Delta_{K}^{+}$ for $K$, can be indexed by Weyl chambers which makes all roots in $\Delta_{K}^{+}$ positive. Given an infinitesimal character $\lambda$ and a Weyl chamber $C$, let $A_{C}(\lambda)$ be the corresponding (L)DS. This is in accordance with the notation of \cite{VZ}. In particular, $A_{C}(\lambda)=A_{\fq}(\lambda_{C})$ for $\fq=\fk\oplus\fp_{C}$, where $\fp_{C}$ is the subalgebra of noncompact roots in $\Delta_{C}^{+}$, the system of positive roots determined by $C$, and $\lambda_{C}$ is the Weyl conjugate of $\lambda$ that is contained in $C$. By the relatively straightforward nature of $K$-multiplicities of such representations, we have the following formulae.
\begin{prop2}\label{Kmultiplicity}
Given $\lambda$ and $C$ as above, we have the following,
\[i_{A_{C}(\lambda)}=\dim_{\bC}(\fp_{-}\cap\fp_{C}),\]
\[V_{A_{C}(\lambda)}=V(\lambda_{C}+2\rho(\fp_{+})),\]the highest weight representation of $K$ with highest weight $\lambda_{C}+2\rho(\fp_{+})$, where $2\rho(\fp_{+})$ is the sum of all roots in $\fp_{+}$. 
\end{prop2}
Therefore, if $\lambda$ lies on some wall of the Weyl chambers, then it can happen that $V_{A_{C}(\lambda)}=V_{A_{C'}(\lambda)}$ for different Weyl chambers $C,C'$. This is the setting we will be interested in. In such cases, we drop the subscript $C$ if there is no issue of confusion.
\begin{thm2}\label{Reducible}
Let $\lambda$ be a nondegenerate analytically integral character. Let $\fC_{\lambda}=\lbrace C\text{ Weyl chamber} \mid \lambda\in C\rbrace$, and $\fP_{\lambda}=\lbrace A_{C}(\lambda)\mid C\in \fC_{\lambda}\rbrace$. Then, there is a parabolic subgroup $Q\subset G(\bR)$ and a discrete series representation $\rho$ of the Levi $M_{Q}$ such that 
\[I(\lambda):=\Ind_{Q}^{G}\rho=\bigoplus_{\pi\in\fP_{\lambda}}\pi,\]where $\Ind_{Q}^{G}\rho$ is the normalized induction. These satisfy the following properties.
\begin{enumerate}
\item $Q$ is a parabolic subgroup which is minimal with respect to the property that the Langlands parameter $\varphi:W_{\bC/\bR}\rar{}^{L}G$ corresponding to the representations in $\fP_{\lambda}$ can be arranged so that $\varphi(W_{\bC/\bR})\subset{}^{L}Q$. Furthermore, if we denote $T_{M_{Q}}$ by a Cartan subgroup of $M_{Q}$, then one can further assume that $\varphi(\bC^{\times})\subset{}^{L}T_{M_{Q}}$ and that $\varphi(W_{\bC/\bR})$ normalizes ${}^{L}T_{M_{Q}}$.
\item $|\fC_{\lambda}|$ is a power of $2$, and $\lambda^{\perp}\subset\fg$ is spanned by a superorthogonal set of real roots.
\item The infinitesimal character of $\rho$ is the restriction of $\lambda$. 
\end{enumerate}
\end{thm2}
\begin{proof}
We freely use the terminology of \cite{Knapp}. The property (3) is clear. The content of \cite[\S14.15]{Knapp} implies that any NLDS appears as a direct summand of such principal series with multiplicity one; this is via repeated application of generalized Schmid identities \cite[Theorem 14.68]{Knapp}. The relation between $Q$ and the Langlands parameter follows from the discussion before \cite[Lemma 1]{LanglandsKZ}. Since the parabolic subgroup appears as Cayley transform of the minimal parabolic in the sense of \cite{KZCorvallis}, and one chooses noncompact simple roots for the Cayley transform, which in fact form a superorthogonal set of roots by \cite[Theorem 14.64]{Knapp}, based only on the infinitesimal character and not the chamber, it follows that \[\lbrace\Ind_{Q}^{G}\wt{\rho}\mid\wt{\rho}\text{ DS, infinitesimal character $=$ the restriction of $\lambda$}\rbrace,\]sees every NLDS with infinitesimal character $\lambda$ at least once as its constituent. The number of constituents is exactly $|W_{G}|/|W_{K}|$, so each such NLDS appears exactly once. Each $\Ind_{Q}^{G}\rho$ has constituents $\pi=A_{C}(\lambda)$ such the representative of $\lambda$ (thought as a Weyl orbit) in $C$ is a fixed character (namely, $C$'s that appear are all adjacent to a single representative of the Weyl orbit of $\lambda$). This is an equivalence relation, so this implies $I(\lambda)$ is precisely consisted of $A_{C}(\lambda)$ where $C$ contains $\lambda$ (regarded as an actual character). As the $R$-group is a direct sum of bunch of $\bZ/2\bZ$ by \cite[\S14.15]{Knapp}, we get (2).
\end{proof}
\begin{exam2}\label{ExampleReducible}We explain how Theorem \ref{Reducible} is realized in some examples. In the figures, red arrows are the compact roots, so NLDS's are those lying on a wall not orthogonal to red arrows.
\begin{enumerate}
\item $\SL_{2}(\bR)$. Let $Q\subset\SL_{2}(\bR)$ be the upper triangular Borel, and let $\rho=\det/|\det|:Q\rar\lbrace\pm1\rbrace$. Then $I_{Q}^{\SL_{2}(\bR)}\rho=D_{0}^{+}\oplus D_{0}^{-}$, the sum of the two NLDS, holomorphic (``weight $1$'') and anti-holomorphic (``weight $-1$'').
\item $\Sp_{4}(\bR)$ (e.g. \cite{Muic}). There are four types of LDS, two of them being holomorphic and anti-holomorphic, respectively, and the other two being large (i.e., maximal Gelfand--Kirillov dimension). We call the one adjacent to the holomorphic chamber \emph{generic} and the other one adjacent to the anti-holomorphic chamber \emph{anti-generic}. Using the notation of \cite{Muic}, a singular infinitesimal character is of one of the forms $(p,0)$, $(0,-p)$ or $(p,-p)$ with $p\in\bN$. The Langlands parameter $\varphi_{\lambda}$ for infinitesimal character $\lambda=(a,b)$ is given by \[\varphi_{\lambda}(re^{i\theta})=\diag(e^{i(a+b)\theta},e^{i(a-b)\theta},e^{-i(a+b)\theta},e^{-i(a-b)\theta}),\quad \varphi_{\lambda}(j)=\begin{pmat}0&(-1)^{a+b}I_{2}\\I_{2}&0\end{pmat}.\]\begin{center}\includegraphics[width=8cm]{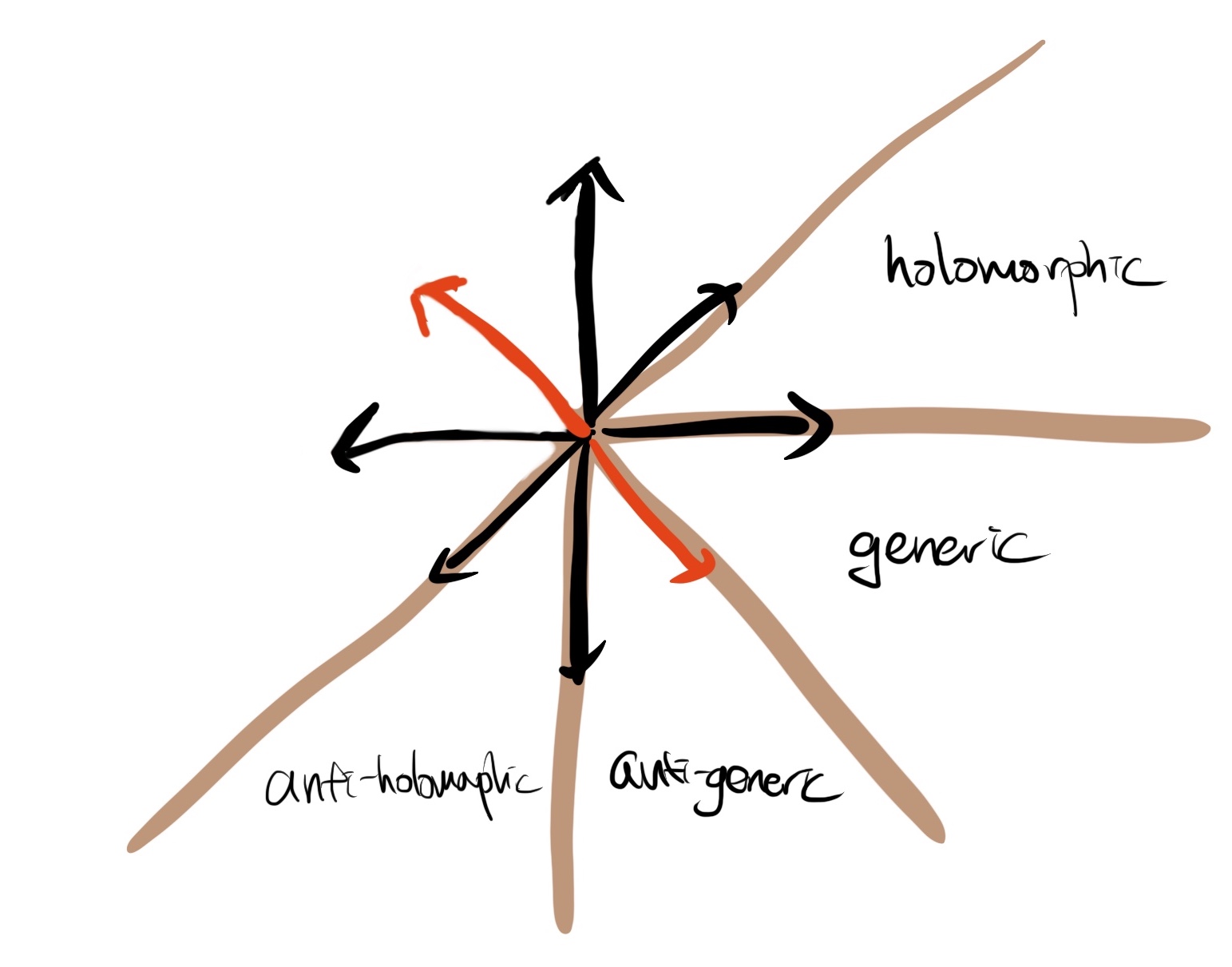}\end{center}\begin{itemize}\item $\lambda=(p,0)$ lies on the wall between the holomorphic chamber and the generic chamber. Then, $Q$ is the so-called \emph{Klingen parabolic}, whose Levi is $\GL_{1}(\bR)\times\SL_{2}(\bR)$, and $\rho$ is the (trivial extension of the) holomorphic DS $D_{p}^{+}$ of $\SL_{2}(\bR)$ of infinitesimal character $p$ (or, equivalently, weight $p+1$).
\item $\lambda=(p,-p)$ lies on the wall between the generic chamber and the anti-generic chamber. Then, $Q$ is the so-called \emph{Siegel parabolic}, whose Levi is $\GL_{2}(\bR)$, and $\rho$ is the DS $D_{2p}$ of $\GL_{2}(\bR)$ with central charcater $2p$ (which is as $\SL_{2}(\bR)$-representation the same as $D_{2p}^{+}\oplus D_{2p}^{-}$).
\item $\lambda=(0,-p)$ lies on the wall between the anti-generic chamber and the anti-holomorphic chamber. This situation is complex-conjugate to the situation of $\lambda=(p,0)$. Thus, $Q$ is again the Klingen parabolic, but $\rho$ is the anti-holomorphic DS of the same infinitesimal character.
\end{itemize} 
\item $\UU(2,1)$ (e.g. \cite{Wallach}, \cite[\S12]{Rogawski}). There are three types of LDS, \emph{holomorphic}, \emph{generic} and \emph{anti-holomorphic}. There are two typse of NLDS's , those lying on the wall between the holomorphic chamber and the generic chamber, and those lying on the wall betweent the generic chamber and the anti-holomorphic chamber.  The Langlands parameters for (L)DS are of the form
\[\varphi(z)=\diag((z/\ov{z})^{a},(z/\ov{z})^{b},(z/\ov{z})^{c}),\quad\varphi(j)=\begin{pmat}&&1\\&-1\\1\end{pmat},\]for $a,b,c\in\bZ$, and the parameter only depends on the unordered set $\lbrace a,b,c\rbrace$. Ordering $a\ge b\ge c$, the two types of NLDS can occur for $a=b>c$ and $a>b=c$. \begin{center}\includegraphics[width=8cm]{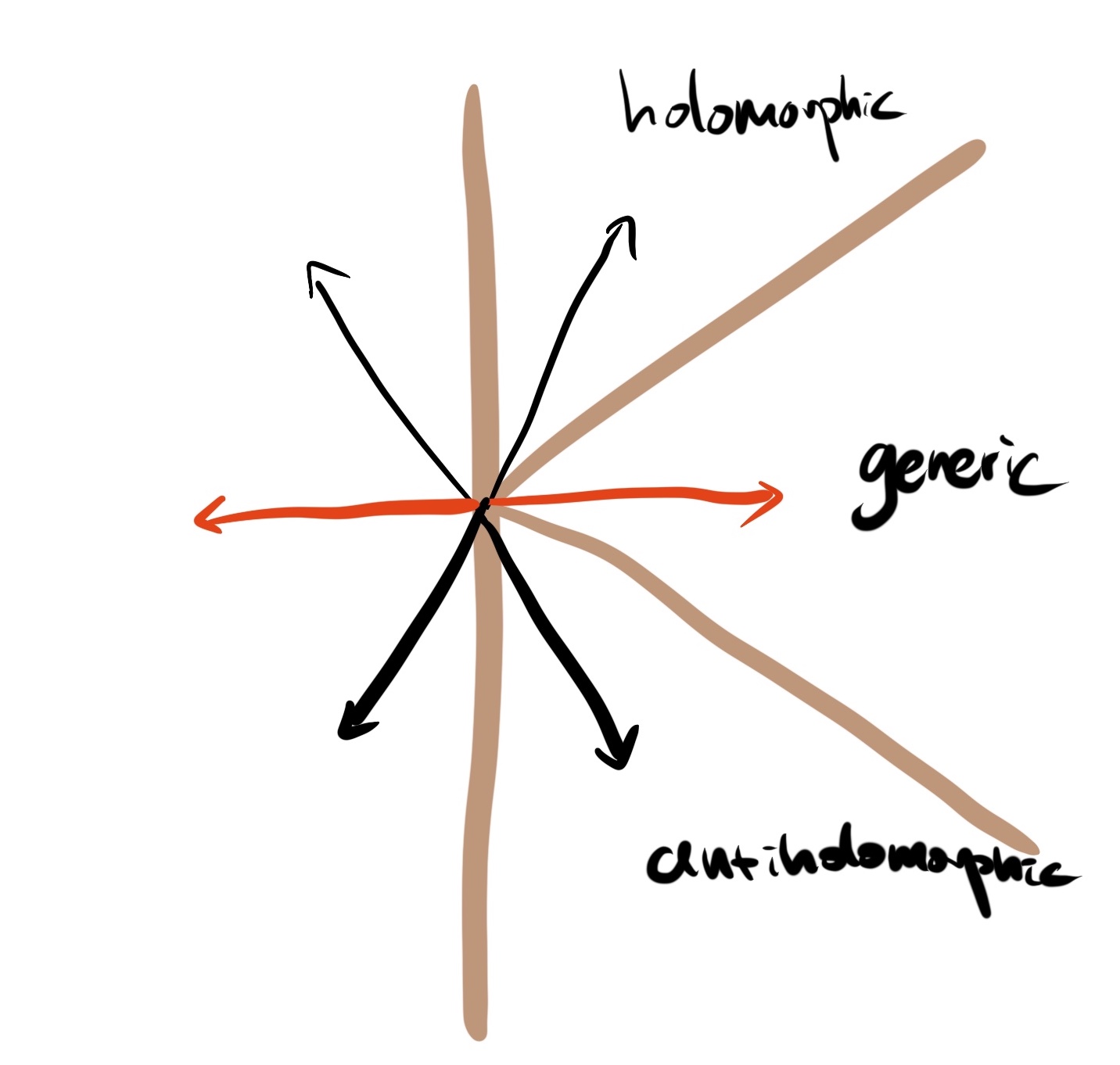}\end{center}

\begin{itemize}
\item $a=b>c$ lies on the wall between the holomorphic chamber and the generic chamber. Then, $Q$ is the upper-triangular Borel (when $\UU(2,1)$ is seen as the unitary group for the diagonal Hermitian matrix such as $\diag(1,1,-1)$), whose Levi is $\bC^{\times}\times S^{1}$, and $\rho$ is the character $\diag(\alpha,\beta,\ov{\alpha}^{-1})\mapsto\alpha^{a}\beta^{b}\ov{\alpha}^{-c}$ for $\alpha\in\bC^{\times}$, $\beta\in S^{1}$.
\item $a>b=c$ lies on the wall between the generic chamber and the anti-holomorphic chamber. The situation is completely symmetric to the previous case.
\end{itemize}
\end{enumerate} 
\end{exam2}
\subsection{Action of the $\Ext$-space}We can now build an archimedean realization of motivic action from abstract nonsense. Recall that for a $(\fp,K)$-module $M$, $H^{i}(\fp,K;M)$ can be regarded as $\Ext_{(\fp,K)}^{i}(\mathbf{1},M)$, the $i$-th $\Ext$ group in the category of $(\fp,K)$-modules, where $\mathbf{1}$ is the trivial module. Thus, there is a natural action
\[\Ext_{(\fp,K)}^{m}(M,N)\times\Ext_{(\fp,K)}^{n}(\mathbf{1},M)\rar\Ext^{m+n}_{(\fp,K)}(\mathbf{1},N).\]
In particular, if $\lambda$ is a nondegenerate singular character as above, with $\pi_{1},\pi_{2}\in\fP_{\lambda}$ with $i_{\pi_{1}}<i_{\pi_{2}}$, there is a natural action
\[\Ext_{(\fp,K)}^{i_{\pi_{2}}-i_{\pi_{1}}}(\pi_{1},\pi_{2})\times H^{i_{\pi_{1}}}(\fp,K;\pi_{1}\otimes V_{A({\lambda})})\rar H^{i_{\pi_{2}}}(\fp,K;\pi_{2}\otimes V_{A_{\lambda}}),\]where $V_{A({\lambda})}$ is the common coefficient for all $\pi\in\fP_{\lambda}$. Thanks to Theorem \ref{Reducible}, the situation can be vastly simplified.
\begin{thm2}\label{Ext}
Let $\lambda$ be as above. \begin{enumerate}\item There is a unique $\pi_{l}\in\fP_{\lambda}$ such that $i_{\pi_{l}}=i_{l}:=\min\lbrace i_{\pi}\mid \pi\in\fP_{\lambda}\rbrace$. \item For each $\pi\in\fP_{\lambda}$, $\Ext^{i_{\pi}-i_{\pi_{l}}}_{(\fp,K)}(\pi_{l},\pi)$ is one-dimensional. \item Let $n_{\lambda}=\log_{2}|\fP_{\lambda}|$, and $I_{\lambda}=\lbrace \pi\in\fP_{\lambda}\mid i_{\pi}=i_{\pi_{l}}+1\rbrace$. Then, $n_{\lambda}=\# I$.

\end{enumerate}
\end{thm2}
\begin{proof}As in the proof of Theorem \ref{Reducible}, \cite[\S14.15]{Knapp} gives us a recipe of how $\lambda$ is constructed in terms of superorthogonal roots $\lbrace\alpha_{1},\cdots,\alpha_{n_{\lambda}}\rbrace$, it follows that the chambers in $\fC_{\lambda}$ are of the following form: there is one chamber $C_{l}$ that is the ``most holomorphic'' chamber (namely, $\dim(C_{l}\cap\fp^{+})$ is the largest) and a set of superorthogonal noncompact roots $\lbrace\alpha_{1},\cdots,\alpha_{n}\rbrace\subset C_{l}$ such that all other chambers $C\in\fC_{\lambda}$ are of the form \[C_{\lbrace i_{1},\cdots,i_{k}\rbrace}:=\left(C_{l}-\lbrace\alpha_{i_{1}},\cdots,\alpha_{i_{k}}\rbrace\right)\cup\lbrace-\alpha_{i_{!}},\cdots,-\alpha_{i_{k}}\rbrace.\]These imply (1) and (3). To show (2), we use the Hochschild--Serre spectral sequence for Lie algebra cohomology, \cite[Proposition 6.1.29]{Vogan}. In our setting, we apply a spectral sequence
\[E_{2}^{p,q}=\Ext_{(\fk,K)}^{p}(X,H^{q}(\fp_{-},Y))\Rar\Ext_{(\fp,K)}^{p+q}(X,Y),\]for $X=\pi_{l}$, $Y=\pi$. Since the category of $(\fk,K)$-modules is semisimple, there is no higher Ext and the spectral sequence trivially degnerates at the $E_{2}$-page:
\[\Ext^{m}_{(\fp,K)}(\pi_{l},\pi)=\Hom_{(\fk,K)}(\pi_{l},H^{m}(\fp_{-},\pi)).\]
On the other hand, the main result of \cite{Williams} says that, if we choose the maximal nilpotent subalgebra $\fb_{-}\subset\fk$ such that $\fn:=\fb_{-}\oplus\fp_{-}$ is a maximal nilpotent subalgebra of $\fg$, then $H^{*}(\fn,\pi)$ can only have weights $W_{K}\lambda+\rho$, where $\rho$ is the half sum of positive roots in $\fg$, where the positivity is defined so that $\fn$ is spanned by the negative roots. Furthermore, each such weight occurs in exactly one cohomological degree. To relate this result to our setting, we use another Hochschild--Serre spectral sequence,
\[E_{2}^{p,q}=H^{p}(\fb_{-},H^{q}(\fp_{-},\pi))\Rar H^{p+q}(\fn,\pi).\]Suppose $V_{\tau}$ be a $K$-type that appears in $H^{m}(\fp_{-},\pi)$, with highest weight $\tau$. Then, by Kostant's theorem \cite[Theorem 5.14]{Kostant}, $H^{p}(\fb_{-},V_{\tau})$ has a nonzero $\mu$-isotypic part if and only if $\mu=w(\tau+\rho_{\fk})+\rho_{\fk}$ for some $w\in W_{K}$ of length $p$, where $\rho_{\fk}$ is the half sum of positive roots in $\fk$, where the positive system is defined so that $\fb_{-}$ is spanned by the negative roots in $\fk$. In particular, $\mu$ determines $p$ and $\tau$, which implies that this spectral sequence also degenerates at the $E_{2}$-page. This implies that the $K$-types of $H^{*}(\fp_{-},\pi)$ are of highest weight $\lambda+\rho_{w\fp_{-}}-\rho_{\fk}$ for $w\in W_{G}/W_{K}$, and these are precisely the minimal $K$-types of $\pi'\in\fP_{\lambda}$.

Let $\tau_{l}$ be the highest weight of the minimal $K$-type of $\pi_{l}$. Then, all other minimal $K$-types of $\pi'\in\fP_{\lambda}$ are of highest weight of the form $\tau_{l}-\sum_{i}^{n_{\lambda}}\epsilon_{i}\alpha_{i}$ for $\epsilon_{i}\in\lbrace0,1\rbrace$. By Blattner's formula for $\pi_{l}$ \cite{Schmid}, all $K$-types of $\pi_{l}$ are of highest weight $\tau_{l}+\delta$, where $\delta$ is a positive linear combination of noncompact positive simple roots, where positivity is defined by $C_{l}$. In particular, this is a cone lying in the direction towards $\alpha_{1},\cdots,\alpha_{n_{\lambda}}$, so the only $K$-type that appears both in $\pi_{l}$ and $H^{m}(\fp_{-},\pi)$ are $\tau_{l}$, whence $\dim\Ext_{(\fp,K)}^{m}(\pi_{l},\pi)\le1$. We have an exact formula for when nonzero Ext group occurs, which turns out to be $i_{\pi}-i_{\pi_{l}}$. Note that the same argument can be used to deduce the following: if $I\subset J\subset\lbrace 1,2,\cdots,n_{\lambda}\rbrace$, then $\Ext^{\#J-\#I}_{(\fp,K)}(\pi_{I},\pi_{J})$ is one-dimensional.\end{proof}
\begin{prop2}\label{ExtHom}

For $I\subsetneq\lbrace1,\cdots,n_{\lambda}\rbrace$ with $i\notin I$, there is a natural identification
\[\Ext_{(\fp,K)}^{1}(\pi_{I},\pi_{I\cup\lbrace i\rbrace})\cong\Hom_{\bC}(\pi_{I}^{\new},\pi_{I\cup\lbrace i\rbrace}^{\new}).\]
\end{prop2}
\begin{proof}

We think the Ext group as the group of actual extensions. In particular, $\Ext_{(\fp,K)}^{1}(\pi_{I},\pi_{I\cup\lbrace\alpha_{i}\rbrace})$ is naturally identified with the one-dimensional vector space of the short exact squences of $(\fp,K)$-modules,
\[0\rar \pi_{I\cup\lbrace\alpha_{i}\rbrace}\rar V\rar \pi_{I}\rar0.\]Such extension is uniquely identified by a scalar $\lambda\in\bC$ in a way as follows. The sequence splits as a short exact sequence of $K$-modules, so the minimal $K$-type $\tau_{I}$ of $\pi_{I}$ has a canonical lift in $V$, also denoted as $\tau_{I}$. Let $v\in\tau_{I}$ be the chosen highest weight vector. As $\alpha_{i}\notin I$, $v$ as a vector in $\pi_{I}$ satisfies $\alpha_{i}\cdot v=0$, where $\alpha_{i}$ is seen as an element of $U(\fp)$. However, if $v$ is seen as a vector in $V$, $\alpha_{i}\cdot v$ is nonzero precisely when the extension is nonsplit, and is sent to the line of highest weight vector of the minimal $K$-type of $\pi_{I\cup\alpha_{i}\rbrace}$. Upon fixing such a highest weight vector $v'$, we have $\alpha_{i}\cdot v=\lambda v'$ for some $\lambda\in\bC$, and this gives an isomorphism $\Ext^{1}_{(\fp,K)}(\pi_{I},\pi_{I\cup\lbrace\alpha_{i}\rbrace})\cong\bC$.
\end{proof}

\begin{rmk2}[(An important new difficulty in the $\delta=0$ case)]\label{Extalgebra}
It is very important to notice that, in our setting, the motivic action depends on the choice of \emph{highest weight vectors of minimal $K$-types} of NLDS's, which was not necessary in the $\delta\ne0$ case of \cite{PV}. This is because our setting involves many different representations, whereas in the setting of \cite{PV} one deals with a single representation. This will give rise to the conjecture on ``generalized complex conjugations'' in Section \ref{GCCSection}.

In the same vein, what naturally acts on the cohomology is not an $\Ext$-algebra, but merely the $\Ext$-group $\Ext^{1}(I(\lambda),\pi_{h})$. However, one can identify this with an exterior algebra as follows. There is a natural identification of $(\fp,K)$-cohomology of $\pi_{I}$ for $I\subset\lbrace1,\cdots,n_{\lambda}\rbrace$,
\begin{align*}H^{*}(\fp,K;\pi_{I}\otimes V_{A(\lambda)})&=\Hom_{K}(\wedge^{*}\fp_{-}, \pi_{I}\otimes V_{A(\lambda)})\\&=\Hom_{K}(\wedge^{*}\fp_{-}\otimes \sigma_{\emptyset},\sigma_{I})\\&=(\wedge^{*}\fp_{+})^{\alpha_{I}},\end{align*}where $\sigma_{J}$ is the minimal $K$-type of $\pi_{J}$, and $\alpha_{I}=\sum_{i\in I}\alpha_{i}$. In particular, this is identified with $(\wedge^{*}\cE_{\lambda})^{\alpha_{I}}$, where $\cE_{\lambda}\subset\fp_{+}$ is the span of $\lbrace\alpha_{1},\cdots,\alpha_{n_{\lambda}}\rbrace$, because these form a superorthogonal set of roots. Thus, one can think of the exterior algebra $\bigwedge^{*}\cE_{\lambda}$ as acting on the $(\fp,K)$-cohomologies of $\pi$'s in $\fP_{\lambda}$. Note that, however, the choice of highest weight vectors of minimal $K$-types underlies everything. 

We will also occasionally use another notation, $\lambda^{\perp}$, for $\cE_{\lambda}$, to be more indicative of its definition.
\end{rmk2}
The action can be connected back to coherent cohomology of Shimura varieties. Consider the $\Pi_{f}$-isotypic part of the coherent cohomology of $X_{G}(\Gamma)$ with coefficient $[V_{A(\lambda)}]$; by Theorem \ref{Jun}, we have a natural isomorphism
\[H^{i_{\Pi_{\infty}}}(X)[\Pi_{f}]\cong H^{i}(\fp,K;\Pi_{\infty}\otimes V_{A(\lambda)})\otimes_{\bC}\Pi_{f}^{\Gamma},\]for each $\Pi_{\infty}\in\fP_{\lambda}$. Thus, using the notation of Theorem \ref{Ext}, we obtain an action
\[\wedge^{*}\cE_{\lambda}\otimes_{\bC}H^{i_{l}}(X)[\Pi_{f}]\rar H^{i_{l}+*}(X)[\Pi_{f}].\]
\subsection{Deligne cohomology as an $\Ext$-space}\label{DeligneExtSection}
We now relate the Ext-space with motivic cohomology. 

By \cite[(2.2.8)]{PV} and the definition of adjoint motive in Conjecture \ref{AdjointMotive}, there is a natural isomorphism $H_{\sD}^{1}((\Ad\Pi)_{\bC},\bR(1))\cong H_{B}((\Ad\Pi)_{\bC},\bR)^{W_{\bC/\bR}}$, and its complexification $H_{B}((\Ad\Pi)_{\bC},\bC)^{W_{\bC/\bR}}$ can be naturally identified with $\wh{\fg}^{\varphi(W_{\bC/\bR})}$, the centralizer of the image of $W_{\bC/\bR}$ under the Archimedean Langlands parameter $\varphi:W_{\bC/\bR}\rar{}^{L}G$. 

We arrange $\varphi$ in a way that is done in Theorem \ref{Reducible}. Then, the centralizer of $\varphi(\bC^{\times})$ is the Levi component $\wh{M_{Q}A_{Q}}$, where $\Lie\wh{M_{Q}A_{Q}}=\wh{\ft}\oplus\wh{\lambda^{\perp}_{\bC}}$. Thus, $\Lie\wh{M_{Q}}$ is a direct sum of $n_{\lambda}$-copies of $\fsl_{2}(\bC)$. Since the centralizer of $\begin{psmat}0&-1\\1&0\end{psmat}$ in $\fsl_{2}(\bC)$ is a one-dimensional torus, $\wh{\fg}^{\varphi(W_{\bC/\bR})}=\bC^{\oplus n_{\lambda}}$ is an $n_{\lambda}$-dimensional abelian Lie algebra, one coming from each $\alpha_{i}$. It is therefore noncanonically isomorphic to $\lambda^{\perp}$. This also shows that $\wh{\fg}^{\varphi(W_{\bC/\bR})}$ gives a Cartan subalgebra $\wh{\fh}=\wh{\ft}\cap\Lie\wh{M_{Q}}$ of $\Lie\wh{M_{Q}}$.

There is on the other hand a natural map $\wh{\fg}^{\varphi(W_{\bC/\bR})}\rar\Ext_{(\fg,K)}^{1}(I(\lambda),\pi_{h})$ whose definition does not depend on any arbitrary choice. Namely, twisting $\sigma$ by a character of $H\subset M_{Q}$ gives a family of deformations of $I(\lambda)=\Ind_{Q}^{G}(\sigma)$, and the same deformation can be applied to extensions: namely, taking $\Ind_{Q}^{G}$ of $\Ext_{\fh}^{*}(\mathbf{1},\mathbf{1})$ gives rise to a map \[\Ext_{\fh}^{*}(\mathbf{1},\mathbf{1})\rar\Ext_{(\fg,K)}^{*}(\Ind_{Q}^{G}(\sigma),\Ind_{Q}^{G}(\sigma)).\]Note that $\Ext_{\fh}^{*}(\mathbf{1},\mathbf{1})=\bigwedge^{*}\wh{\fh}$, so we have a natural map
\[\wedge^{*}\wh{\fh}=\Ext_{\fh}^{*}(\mathbf{1},\mathbf{1})\rar\Ext_{(\fg,K)}^{*}(I(\lambda),I(\lambda))=\bigoplus_{\pi\in\fP_{\lambda}}\Ext_{(\fg,K)}^{*}(I(\lambda),\pi)\thrar\Ext_{(\fg,K)}^{*}(I(\lambda),\pi_{h}),\]where the last map is the projection map. This isomorphism is natural, so this defines a \emph{degree-descending} action of $\bigwedge^{*}H_{\sD}^{1}(\Ad\Pi,\bC(1))$ on $H^{i_{h}-*}(X)[\Pi_{f}]$. Dually, this defines a \emph{degree-ascending} action of $\bigwedge^{*}H_{\sD}^{1}(\Ad\Pi,\bC(1))^{*}$ on $H^{i_{l}+*}(X)[\Pi_{f}]$.
\bibliographystyle{alpha}

\end{document}